\documentclass[12pt]{article}
\usepackage{amssymb,amsmath,theorem}
\usepackage[all]{xy}
\usepackage{lscape}

\title{\bf The cohomology ring of free loop spaces}
\date{}
\author{Luc Menichi}

\newtheorem{theor}{Theorem}[section]
\newtheorem{cor}[theor]{Corollary}
\newtheorem{proposition}[theor]{Proposition}
\newtheorem{lem}[theor]{Lemma}
{\theorembodyfont{\rmfamily}\newtheorem{defin}[theor]{Definition}
\newtheorem{rem}[theor]{Remark}
\newtheorem{ex}[theor]{Example}
\newtheorem{propriete}[theor]{Property}
}

\def\QED{\hskip.1em\hfill\vrule width0pt height2pt\nobreak\hfill
   \kern3pt\lower1.8pt\vbox{\hrule\hbox
   {\vrule\kern1pt\vbox{\kern1.7pt\hbox{$\scriptstyle
   QED$}\kern0.2pt}\kern1pt\vrule}\hrule}}

\newcommand{\pf}{\noindent{\bf Proof.}\ }

\newcommand{\T}{\mathrm{T}}
\newcommand{\TA}{\mathrm{T\!A}}
\newcommand{\TC}{\mathrm{T\!C}}

\newcommand{\B}{\mathrm{B}}
\newcommand{\C}{\mathrm{C}}

\def\build#1_#2^#3{\mathrel{\mathop{\kern0pt#1}\limits_{#2}^{#3}}}

\begin{document}
\maketitle
\begin{abstract}
Let $X$ be a simply connected space and $\Bbbk$ a commutative ring.
Goodwillie, Burghelea and Fiedorowiscz proved that the Hochschild
cohomology of the singular chains on the pointed loop space
$HH^{*}S_*(\Omega X)$ is isomorphic to the free loop space
cohomology $H^{*}(X^{S^{1}})$.
We proved that this isomorphism is compatible with both
the cup product on $HH^{*}S_*(\Omega X)$ and on $H^{*}(X^{S^{1}})$.
In particular, we explicit the algebra $H^{*}(X^{S^{1}})$
when $X$ is a suspended space, a complex projective space
or a finite CW-complex of dimension $p$ such
that $\frac {1}{(p-1)!}\in {\Bbbk}$.
\end{abstract}

Research supported by the University of Toronto (NSERC grants RGPIN 8047-98 and OGP000 7885).
\newpage
\section{Introduction}
Let $X$ be a simply connected CW-complex.
We consider the free loop space cohomology algebra of $X$,
$H^{*}(X^{S^{1}})$, with coefficients in any arbitrary commutative
ring ${\Bbbk}$.

When ${\Bbbk}$ is a field, the algebra $H^{*}(X^{S^{1}})$ can
sometimes
be computed via spectral sequences.
For example, Kuribayashi and Yamaguchi~\cite{Kuribayashi},
by solving extensions problems by applications of the Steenrod
operations, were able to compute via the Eilenberg-Moore spectral
sequence, the algebra $H^{*}(X^{S^{1}})$ for some simple spaces.

If $A$ is a commutative algebra, the Hochschild homology of $A$,
$HH_*(A)$ can be endowed with the shuffle product.
When ${\Bbbk}=\mathbb Q$, Sullivan, inventing Rational Homotopy,
constructed a commutative algebra $A_{PL}(X)$, called the
{\it polynomial differential forms}~\cite[\S 10]{rational homotopy},
and proved with Vigu\'e-Poirrier~\cite{Sullivan et Micheline}
that there is a natural isomorphism of graded algebras
$$H^{*}(X^{S^{1}})\cong HH_*(A_{PL}(X)).$$
The theory of minimal Sullivan models gives a technic to
compute the algebra $HH_*(A_{PL}(X))$ and so $H^{*}(X^{S^{1}})$.

Denote by $HH^{*}(A)$ the Hochschild cohomology of an algebra $A$
with coefficients in $A^{\vee}=\mbox{Hom}(A,{\Bbbk})$~\cite[1.5.5]{Loday}.
If $A$ has a diagonal $\Delta:A\rightarrow A\otimes A$,
then $HH^{*}(A)$ is equipped with a cup product.
The normalized singular chains on the pointed loop space
$S_*(\Omega X)$ is a differential graded Hopf algebra.
So the Hochschild cohomology $HH^{*}(S_*(\Omega X))$ is naturally 
a graded algebra.
In 1985, Goodwillie~\cite{Goodwillie}, Burghelea and
Fiedorowiscz~\cite{Burghelea-Fiedorowiscz} proved that there
is an isomorphism of graded modules
$$HH^{*}((S_*(\Omega X))\cong H^{*}(X^{S^{1}}).$$

The main result of this paper is that this isomorphism
of graded modules is, in fact, an isomorphism of graded algebras
(Theorem~\ref{Theoreme de Goodwillie}).
This result gives a general method to compute the algebra
$H^{*}(X^{S^{1}})$
over any commutative ring:
\begin{itemize}
\item Determine the Adams-Hilton model \cite{modeles d'adams-hilton}
of $X$, ${\cal A}(X)$, using its cellular decomposition.
\item Compute the structure of Hopf algebra up to homotopy model on
  ${\cal A}(X)$
(Definition~\ref{definition HAH model}), knowing the structure of
graded Hopf algebra of $H_*(\Omega X)$ when ${\Bbbk}$ is a field.
\end{itemize}
Note that, if $X$ is an H-space, our approach has no interest,
since simply $$X^{S^{1}}\thickapprox\Omega X\times X.$$
Now, Theorem~\ref{cohomologie lacets libres et HAH}
allows us to replace, in the Hochschild cohomology, the
differential graded Hopf algebra $S_*(\Omega X)$
by the Adams-Hilton model of $X$, ${\cal A}(X)$,
equipped with its structure of Hopf algebras up to homotopy.

We face now (section~\ref{smaller coalgebra}) a completely algebraic problem:
How to compute the Hochschild cohomology
on a Hopf algebra up to homotopy
whose underlying algebra is a tensor algebra.

In section~\ref{homologie de Hochschild et Cobar},
we show that in a simple case, this Hochschild cohomology reduces
to the Hochschild homology of a commutative algebra.

In section~\ref{free loop space on a suspension},
we investigated the algebra structure of the free loop space
cohomology on any suspension $\Sigma X$.
If $H_*(X)$ is ${\Bbbk}$-free of finite type, this algebra
$H^{*}\left((\Sigma X)^{S^{1}})\right)$ is the Hochschild
homology of $H^{*}(\Sigma X)$, $HH_*(H^{*}(\Sigma X))$,
equipped with a product completely determined by the cohomology algebra
of the desuspended space, $H^{*}(X)$.
Recall from~\cite{Lambe et Lofwall} and~\cite{Parhizgar} that,
even when ${\Bbbk}=\mathbb{Q}$ and $\Sigma X$ is a wedge of spheres,
the cohomology algebra $H^{*}\left((\Sigma X)^{S^{1}})\right)$
is particularly difficult to explicit in terms of
generators and relations.

In section~\ref{exemples Hochschild homology d'une algebre
  commutative},
we prove that the free loop space cohomology on the complex projective
space $\mathbb{CP}^{n}$,
  $H^{*}\left((\mathbb{CP}^{n})^{S^{1}}\right)$,
is isomorphic as graded algebras to the Hochschild homology
$HH_*(H^{*}(\mathbb{CP}^{n}))$ and compute it.
Suppose that $X$ is a finite CW-complex of dimension $p$ such
that $\frac {1}{(p-1)!}\in {\Bbbk}$.
Anick~\cite[dualize Proposition 8.7(a)]{anick},
extending Sullivan's result, constructed a commutative algebra
${\mathcal C}^{*}(\mathbf{L}(X))$ weakly equivalent as algebras
(in the sense of~\cite[page 832]{dga in topology}) 
to the singular cochains on $X$, $S^{*}(X)$.
We extend Sullivan and Vigu\'e-Poirrier result in this new context
(Theorem~\ref{domaine d'Anick}).

We would like to mention that Ndombol and Thomas~\cite{bitjong et thomas}
have also found, when ${\Bbbk}$ is a field, a general method
to compute the free loop space cohomology algebra of a simply connected
space.
We thank S. Halperin, J.-C. Thomas and M. Vigu\'e for their constant
support.
The main results of this paper were exposed in September 1999, at the GDR
Topologie alg\'ebrique meeting in Paris Nord.
\section{Algebraic preliminaries and notation}\label{preliminaries}
We work over a commutative ring ${\Bbbk}$.
We denote by $ _p{\Bbbk}$ and $\frac{\Bbbk}{p{\Bbbk}}$ respectively
the kernel and cokernel of the multiplication by $p$
in ${\Bbbk}$.

DGA stands for differential graded algebra, DGC for differential
graded coalgebra,
DGH for differential graded Hopf algebra and CDGA for commutative DGA.
The denomination ``chain'' will be restricted to objects with a non-negative lower degree
and ``cochain'' to those with a non-negative upper degree.

The degree of an element $x$ is denoted $\vert x\vert$.
The {\it suspension} of a graded module $V$ is the graded module $sV$ such that
$(sV)_{i+1}=V_i$.
Let $C$ be an augmented complex. The kernel of the augmentation is denoted
$\overline{C}$.

The exterior algebra on an element $v$ is denoted $Ev$.
The {\it free divided powers algebra} on an element $v$, denoted $\Gamma v$,
is
\begin{itemize}
\item the free graded algebra generated by
$\gamma^{i}(v),\:i\in \mathbb{N}^{*}$, 
divided by the relations
$\displaystyle\gamma^{i}(v)\gamma^{j}(v)=\frac{(i+j)!}{i!j!}\gamma^{i+j}(v)$,
if $\vert v\vert$ is even,
\item and is just $Ev$ when $\vert v\vert$ is odd.
\end{itemize}
The {\it tensor algebra} on a graded module $V$ is denoted $\TA V$.
The {\it tensor coalgebra} is denoted $\TC V$.
Their common underlying module is simply denoted $\T V$.
Given a conilpotent coalgebra $C$ then any morphism
$\varphi:\overline{C}\rightarrow V$ lifts uniquely to a unique
morphism
$\Psi:C\rightarrow\TC V$ of coaugmented coalgebras.
The formula for $\Psi$ is given by
\addtocounter{theor}{1}
\begin{equation}\label{formule du relevement coalgebre colibre}
\Psi(c)=\sum_{i=1}^{+\infty}
\varphi ^{\otimes i}\circ\overline{\Delta}^{\otimes i-1}_C(c),c\in\overline{C}
\end{equation}
where $\Delta^{i-1}_C:\overline{C}\rightarrow \overline{C}^{\otimes
i}$
is the iterated reduced diagonal of $C$.

Let $A$ be an augmented DGA.
Denote $d_1$ be the differential of the complex $A\otimes \T(s\overline{A})\otimes A$
obtained by tensorization.
We denote the tensor product of the elements $a\in A$, $sa_1\in s\overline{A}$,
\ldots , $sa_k\in s\overline{A}$ and $b\in A$ by $a[sa_1|\cdots|sa_k]b$.
Let $d_2$ be the differential on the graded module
$A\otimes \T(s\overline{A})\otimes A$ defined by:
\begin{align*}
d_2a[sa_1|\cdots|sa_k]b=&(-1)^{\vert a\vert } aa_1[sa_2|\cdots|sa_k]b\\
&+\sum _{i=1}^{k-1} (-1)^{\varepsilon_i}{a[sa_1|\cdots|sa_ia_{i+1}|\cdots|sa_k]b}\\
&-(-1)^{\varepsilon_{k-1}} a[sa_1|\cdots|sa_{k-1}]a_kb;
\end{align*}
Here $\varepsilon_i=\vert a\vert +\vert sa_1\vert+\cdots +\vert sa_i\vert$.

The {\it bar resolution of $A$}, denoted
$\B(A;A;A)$, is the $(A,A)$-bimodule $(A\otimes
\T(s\overline{A})\otimes A,d_1+d_2)$.
The {\it (reduced) bar construction on $A$}, denoted
$\B(A)$, is the coaugmented DGC $(\TC s\overline{A},d_1+d_2)$
whose underlying complex $(\T s\overline{A},d_1+d_2)$ coincides with
${\Bbbk}\otimes_A \B(A;A;A)\otimes_A {\Bbbk}$ \cite[\S 4]{dga in topology}.
The {\it cyclic bar construction} or {\it Hochschild complex} is the complex
$A\otimes_{A\otimes A^{op}} \B(A;A;A)$ denoted $\C(A)$.
Explicitly $\C(A)$ is the complex ($A\otimes \T(s\overline{A}),d_1+d_2)$
with $d_1$ obtained by tensorization and
\begin{align*}
d_2a[sa_1|\cdots|sa_k]=&(-1)^{\vert a\vert } aa_1[sa_2|\cdots|sa_k]\\
&+\sum _{i=1}^{k-1} (-1)^{\varepsilon_i}{a[sa_1|\cdots|sa_ia_{i+1}|\cdots|sa_k]}\\
&-(-1)^{\vert sa_k\vert\varepsilon_{k-1}} a_ka[sa_1|\cdots|sa_{k-1}];
\end{align*}
The {\it Hochschild homology} is the homology of the cyclic bar construction:
$$HH_*(A):=H_*(\C(A)).$$
The {\it Hochschild cohomology} is the graded module
$$HH^{*}(A):=H^{*}(\mbox{Hom}_{(A,A)} (\B(A;A;A),A^{\vee}))=H^{*}(\C(A)^{\vee})$$
where $A^{\vee}$ is considered as an $(A,A)$-bimodule.

Let $A$ and $B$ be two augmented DGA's,
Then we have an Alexander-Whitney morphism of $(A\otimes B, A\otimes B)$-bimodules
$$AW:\B(A\otimes B;A \otimes B;A\otimes B)\rightarrow \B(A;A;A)\otimes \B(B;B;B)$$
where the image of a typical element
$p\otimes q [s(a_1 \otimes b_1)|\cdots |s(a_k \otimes b_k)]m\otimes n$ is
$$\sum _{i=0}^{k}(-1)^{\zeta_i}p[sa_1|\cdots |sa_i]a_{i+1}\cdots a_km
\otimes qb_1\cdots b_i[sb_{i+1}|\cdots |sb_k]n.$$
Here~\cite[3.7]{Cohomology of restricted Lie algebras}
\begin{multline*}
\zeta_i=\sum_{j=1}^{k}\left (\vert q\vert+\sum_{l=1}^{j-1}\vert b_l\vert\right )
\vert a_j\vert+\left (\vert q\vert+\sum_{j=1}^{k}\vert b_j\vert\right )\vert m\vert\\
+\sum_{j=i+1}^{k}(j-i)\vert a_j\vert+(k-i)\vert m\vert+\vert i\vert\vert q\vert
+\sum_{j=1}^{i-1}(i-j)\vert b_j\vert.
\end{multline*}
$AW$ is natural and associative exactly.
It is also commutative up to a homotopy of $(A\otimes B, A\otimes B)$-bimodules.
So we get an Alexander-Whitney map for the cyclic bar construction
$$AW:\C(A\otimes B)\rightarrow \C(A)\otimes \C(B).$$
Consider an augmented DGA $K$ equipped with a morphism of augmented DGA's
$\Delta:K\rightarrow K\otimes K$.
Then the composite
$$\Delta:\C(K)\buildrel{\C(\Delta)}\over\longrightarrow \C(K\otimes K)
\buildrel{AW}\over\longrightarrow \C(K)\otimes \C(K)$$
is a morphism of augmented complex. Therefore $HH^{*}(K)$ has a
product. This the cup product of Cartan and Eilenberg \cite[XI.6]{Cartan-Eilenberg}.
In particular,
if $K$ is a DGH then $\C(K)$ is a DGC and $HH^{*}(K)$ is a graded algebra.

\section{From the chains on the based loops to the chains on the free loops}
The object of this section is to prove the following theorem linking the
chains on the based loops of a space to the chains on its free loops.
\begin{theor}(Compare \cite{Goodwillie} and \cite{Burghelea-Fiedorowiscz})\label{Theoreme de Goodwillie}
Let $X$ be a path connected pointed space.
Then there is a natural DGC quasi-isomorphism
$$\C(S_*(\Omega X))\buildrel{\simeq}\over\rightarrow S_*(X^{S^{1}}).$$
In particular, $HH^{*}S_*(\Omega X)\cong H^{*}(X^{S^{1}})$ as graded algebras.
\end{theor}
Goodwillie~\cite{Goodwillie}, Burghelea and Fiedorowiscz~\cite{Burghelea-Fiedorowiscz}
proved the isomorphism
$HH^{*}S_*(\Omega X)\cong H^{*}(X^{S^{1}})$ as graded modules only.
To obtain our theorem, we will follow their proofs.
We introduce first some terminology about simplicial objects.

Let $\mathcal{C}$ be a category.
A {\it simplicial} $\mathcal {C}$-object $X$ is a non-negative graded object together with morphisms
$d_i:X_n\rightarrow X_{n-1}$ and $s_i:X_n\rightarrow X_{n+1}$, $0\leq i\leq n$ satisfying some well-known
relations~\cite[VIII.5.2]{homology}.
A {\it cosimplicial} $\mathcal {C}$-object is a non-negative graded object together with morphisms
$\delta_i:X^{n-1}\rightarrow X^{n}$ and $\sigma_i:X^{n+1}\rightarrow X^{n}$, $0\leq i\leq n$ satisfying the opposite
relations~\cite[X.2.1(i)]{Bousfield-Kan}.
If  $\mathcal {C}$ is a category equipped with a tensor product $\otimes$
(more precisely a monoidal category~\cite[VII.1]{working mathematician})
then the tensor product of two simplicial $\mathcal {C}$-objects
$X=(X_n, d_i, s_i)$ and $Y=(Y_n, d_i, s_i)$ is the simplicial $\mathcal {C}$-object
$X\otimes Y=(X_n\otimes Y_n, d_i\otimes d_i, s_i\otimes s_i)$.

Consider $\mathcal {C}$ to be the category of complexes.
To any simplicial $\mathcal {C}$-object (i.e. simplicial complex) $X$, we can associate a complex in the category $\mathcal {C}$
(i.e. a complex of complexes) denote $K_{N}(X)$ known as the normalized chain complex of
$X$~\cite[VIII.6 for the category of modules]{homology}.
Consider two simplicial complexes $A$ and $B$. We have an Alexander-Whitney morphism of
complexes of complexes~\cite[VIII.8.6]{homology} $AW:K_N(A\otimes B)\rightarrow K_N(A)\otimes K_N(B)$.
Every complex of complexes can be condensated~\cite[X.9.1]{homology} into a single complex.
So by composing the functor $K_N$ and the condensation functor, we
have a functor, called the {\it realization} and denoted
$\vert\phantom{X}\vert$, from the category
of simplicial complexes to the category of complexes,
equipped with an Alexander-Whitney morphism of
complexes $AW:\vert A\otimes B\vert \rightarrow \vert A\vert\otimes
\vert B\vert$ for any simplicial complexes $A$ and $B$.
In particular, $\vert\phantom{X}\vert$ induces a functor from the category of simplicial DGC's to the category of DGC's
(Recall that a simplicial DGC can be defined either as a simplicial object in the category of DGC's
or as a coalgebra in the category of simplicial complexes.).

Given any two topological spaces $X$ and $Y$, the caligraphic
notations
$$\mathcal{AW}:S_*(X\times Y)\rightarrow S_*(X)\otimes S_*(Y)
\mbox{ and }\mathcal{EZ}:S_*(X)\otimes S_*(Y)\rightarrow S_*(X\times
Y)$$
are reserved to the standart normalized Alexander-Whitney map and to
the
standart normalized Eilenberg-Zilber map concerning singular chains
\cite[VI.12.27-8]{Dold}.
\begin{ex}\label{Bar construction cyclique simpliciale}
The cyclic bar construction for differential graded algebras.
Let $A$ be a DGA. Then there is a simplicial complex $\Gamma A$ defined by
$\Gamma_nA=A\otimes\cdots\otimes A=A^{\otimes n+1}$,
\begin{align*}
d_0a[a_1|\cdots|a_n]&=aa_1[a_2|\cdots|a_n],\\
d_ia[a_1|\cdots|a_n]&=a[a_1|\cdots|a_ia_{i+1}|\cdots|a_n]\mbox{ for } 1\leq i\leq n-1,\\
d_na[a_1|\cdots|a_n]&=a_na[a_1|\cdots|a_{n-1}],\\
s_ia[a_1|\cdots|a_n]&=a[a_1|\cdots|a_i|1|a_{i+1}|\cdots|a_n]\mbox{ for }0\leq i\leq n.
\end{align*}
The complex $\vert\Gamma  A\vert$ is exactly (signs included) $\C(A)$ the cyclic bar construction of $A$.
If $K$ is a DGH then $\Gamma K$ with the diagonal
$$\Gamma K\xrightarrow{\Gamma(\Delta_K)}\Gamma (K\otimes K)\cong\Gamma K\otimes\Gamma K$$
is a simplicial DGC and $\vert\Gamma K\vert$ is the DGC $\C(K)$ denoted in section \ref{preliminaries}.
\end{ex}
\begin{ex}
Let $G$ be a topological monoid. The cyclic bar construction of $G$ \cite[7.3.10]{Loday} is the simplicial
space $\Gamma G$ defined by $\Gamma_n G=G\times\cdots\times G=G^{n+1}$
and with the same formulas for
$d_i$ and $s_i$ as in the cyclic bar construction for DGA's.
Since the normalized singular chain functor $S_*$ is a functor from topological spaces to DGC's,
$S_*(\Gamma G)$ is a simplicial DGC. Therefore $\vert S_*(\Gamma
G)\vert$ is a DGC.
\end{ex}
The following Lemma compares the DGC's given by the previous two examples.
\begin{lem}(Compare~\cite[V.1.2]{Goodwillie})\label{echange cyclic et chaines singulieres}
Let $G$ be a topological monoid. Then there is a natural DGC quasi-isomorphism
$\vert\Gamma S_*(G)\vert\buildrel{\simeq}\over\rightarrow
\vert S_*(\Gamma G)\vert$.
\end{lem}
\pf The Eilenberg-Zilber map $\mathcal{EZ}:S_*(G)^{\otimes n+1}\rightarrow S_*(G^{n+1})$ is a DGC quasi-isomorphism
and therefore defines a morphism of simplicial DGC's $\Gamma S_*(G)\rightarrow S_*(\Gamma G)$.
So, applying the functor $\vert\phantom{X}\vert$,
we get a DGC quasi-isomorphism.\QED

Let $\Delta^{n}$ be the standart geometric simplex of dimension $n$.
Let $\delta_i:\Delta^{n-1}\rightarrow\Delta^{n}$ and
$\sigma_i:\Delta^{n+1}\rightarrow\Delta^{n}$ be the $i$-th face inclusion
and the $i$-th degeneracy of $\Delta^{n}$.
Then $\mathbb{\Delta}=(\Delta^{n},\delta_i,\sigma_i)$
is a cosimplicial space~\cite[X.2.2(i)]{Bousfield-Kan}.
The geometric {\it realization} \cite[11.1]{Geometry of iterated
  loop space} of a simplicial space $X$ is defined as
$$\vert X\vert=\left( \coprod_{n\in\mathbb{N}}X_n\times\Delta^{n} \right)\bigg /\sim$$
where $\sim$ is the equivalence relation generated by
$$\phantom{\mbox{and}\quad}(d_ix,y)\sim (x,\delta_iy),\;x\in X_n,\;y\in \Delta^{n-1}$$
$$\mbox{and}\quad (s_ix,y)\sim (x,\sigma_iy),\;x\in X_n,\;y\in
\Delta^{n+1}.$$
Recall that $\vert S_*(X)\vert$ is a DGC whose diagonal is the
composite
$$\vert S_*(X)\vert \buildrel{\vert S_*(\Delta)\vert}\over\rightarrow
\vert S_*(X\times X)\vert
\buildrel{\vert\mathcal{AW}\vert}\over\rightarrow
\vert S_*(X)\otimes S_*(X)\vert
\buildrel{AW}\over\rightarrow
\vert S_*(X)\vert\otimes\vert S_*(X)\vert$$
\begin{lem}(Compare \cite[Theorem 4.1]{Bott} and \cite[17(a)]{rational homotopy})\label{realisation}
Let $X$ be a simplicial space, good in the sense of \cite[A.4]{Segal}.
Then there is a natural DGC quasi-isomorphism
$f:\vert S_*(X)\vert\buildrel{\simeq}\over\rightarrow S_*(\vert X\vert)$.
\end{lem}
\pf Let $\pi_n:X_n\times\Delta^{n}\twoheadrightarrow\vert X\vert$ be the quotient map.
The morphism $f$ is defined as the composite
$$S_i(X_n)\xrightarrow{id_{S_i(X_n)}\otimes\kappa_n}
S_i(X_n)\otimes S_n(\Delta^{n})\buildrel{\mathcal{EZ}}\over\rightarrow S_{i+n}(X_n\times\Delta^{n})
\xrightarrow{S_{i+n}(\pi_n)}S_{n+i}(\vert X\vert)$$
where $\kappa_n\in S_n(\Delta^{n})$ is the singular simplex $id_{\Delta^{n}}$.
We just have to prove that $f$ is a DGC morphism.
The diagonal map of $\vert X\vert$ is equal to the composite
$$\vert X\vert\xrightarrow{\vert\Delta_X\vert}\vert X\times X\vert
\xrightarrow{(\vert proj_1\vert,\vert proj_2\vert)}\vert X\vert\times\vert X\vert$$
where $\vert\Delta_X\vert$ is the simplicial diagonal of $X$ and
where $proj_1$ and $proj_2$ are the simplicial projections on each factors.
So by naturality of $f$, it suffices to show that $f$ well behaves with products of
simplicial spaces.

Let $X$ and $Y$ be two simplicial spaces. Then $S_*(X)$ and $S_*(Y)$
are two simplicial complexes. So we have an Alexander-Whitney map
$$\vert S_*(X)\otimes S_*(Y)\vert\buildrel{AW}\over\rightarrow
\vert S_*(X)\vert\otimes\vert S_*(Y)\vert.$$
Its formula is given by
$$\left[\sum_{p+q=n}S_*(\tilde{d}^{q})\otimes S_*(d_0^{p})\right]:S_*(X_n)\otimes S_*(Y_n)\longrightarrow
\bigoplus_{p+q=n}S_*(X_p)\otimes S_*(Y_q)$$
where $\tilde{d}^{q}:X_n\rightarrow X_p$ is the composite $d_{p+1}\circ\cdots\circ d_n$
($\tilde{d}$ denotes the ``last'' face operator) and $d_0^{p}:Y_n\rightarrow Y_q$ is
the iterated composite of $d_0$.
In the diagram page~\pageref{coalgebre et realisation}, there was no space left for sums $\sum$ and direct sums $\bigoplus$.
So we use the indices $p$ and $q$ with the convention $p+q=n$ and the indices
$j$ and $k$ with the conventions that $j+k=i$.
We use also the maps
\[\begin{array}{ll}
\tilde{\delta}^{q}=\delta_n\circ\cdots\circ\delta_{p+1}:\Delta^{p}\rightarrow\Delta^{n},
&\delta^{p}_0:\Delta^{q}\rightarrow\Delta^{n},\\
\tilde{\sigma}^{q}_0=\sigma_p\circ\cdots\circ\sigma_{n-1}:\Delta^{n}\rightarrow\Delta^{p}
\quad\mbox{and}
&\sigma^{p}_0:\Delta^{n}\rightarrow\Delta^{q}.
\end{array}
\]
Both the interchange of factors of a tensor product of modules and of a product of spaces
are denoted by $\tau$.

Consider the diagram page~\pageref{coalgebre et realisation}.
\begin{landscape}
\voffset -3cm
\label{coalgebre et realisation}
\xymatrix@C=3pc@R=5pc{
S_i(X_n\times Y_n) \ar[r]^-{id\otimes\kappa_n} \ar[dd]_{\mathcal{AW}}
&S_i(X_n\times Y_n)\otimes S_n(\Delta^{n}) \ar[r]^{\mathcal{EZ}} \ar[d]|{id\otimes S_n(\Delta)}
\ar  @{} [ddl] _*+[o][F-]{1}
\ar @{} [dr] |*+[o][F-]{3}
&S_{n+i}(X_n\times Y_n\times\Delta^{n}) \ar[r]^{S_{n+i}(\pi_n)}
\ar[d]|{S_{n+i}(id\times\Delta)}
\ar @{} [ddr] |*+[o][F-]{6}
&S_{n+i}(\vert X\times Y\vert )
\ar[dd]|{S_{n+i}(\vert proj_1\vert,\vert proj_2\vert)}\\
&S_i(X_n\times Y_n)\otimes S_n(\Delta^{n}\times\Delta^{n})
\ar[d]|{\mathcal{AW}\otimes \mathcal{AW}} \ar[r]^{\mathcal{EZ}}
\ar @{} [dr] |*+[o][F-]{4}
&S_{n+i}(X_n\times Y_n\times\Delta^{n}\times\Delta^{n})
\ar[d]|{S_{n+i}(id\times\tau\times id)}\\
S_{j}(X_n)\otimes S_{k}(Y_n)
\ar@/^/[r]^-{id\otimes \mathcal{AW}(\kappa_n,\kappa_n)}
\ar[d]|{ S_{j}(\tilde{d}^{q})\otimes S_{k}(d_0^{p})}
\ar @{} [dr] |*+[o][F-]{2}
&S_j(X_n)\otimes S_{k}(Y_n)\otimes S_{p}(\Delta^{n})\otimes S_{q}(\Delta^{n})
\ar[d]|{S_{j}(\tilde{d}^{q})\otimes S_{k}(d_0^{p})\otimes S_{p}(\tilde{\sigma}^{q})\otimes S_{q}(\sigma^{p}_0)}
\ar[dr]^{\mathcal{EZ}\otimes \mathcal{EZ}\circ (id\otimes\tau\otimes id)}
\ar @{} [drd] |*+[o][F-]{5}
&S_{n+i}(X_n\times\Delta^{n}\times Y_n\times\Delta^{n})
\ar[r]^-{S_{n+i}(\pi_n\times\pi_n)} \ar[d]^{\mathcal{AW}}
\ar @{} [dr] |*+[o][F-]{7}
&S_{n+i}(\vert X\vert\times\vert Y\vert)
\ar[d]^{\mathcal{AW}}\\
S_{j}(X_{p})\otimes S_{k}(Y_{q})
\ar[r]_-{id\otimes\kappa_p\otimes\kappa_q}
&S_{j}(X_{p})\otimes S_{k}(Y_{q})\otimes S_{p}(\Delta^{p})\otimes S_{q}(\Delta^{q})
\ar[dr]_{\mathcal{EZ}\otimes \mathcal{EZ}\circ (id\otimes\tau\otimes id)}
&S_{p+j}(X_n\times\Delta^{n})\otimes S_{q+k}(Y_n\times\Delta^{n})
\ar@/^/[r]^-{S_{p+j}(\pi_n)\otimes S_{q+k}(\pi_n)}
\ar[d]|{S_{p+j}(\tilde{d}^{q}\times\tilde{\sigma}^{q})\otimes S_{q+k}(d_0^{p}\times\sigma^{p}_0)}
&S_{p+j}(\vert X\vert )\otimes S_{q+k}(\vert Y\vert)
\ar @{} [dl] |*+[o][F-]{8}\\
&& S_{p+j}(X_{p}\times\Delta^{p})\otimes S_{q+k}(Y_{q}\times\Delta^{q})
\ar@/_2pc/[ur]_-{S_{p+j}(\pi_{p})\otimes S_{q+k}(\pi_{q})}
}
\end{landscape}
\voffset 0cm

Let's check the commutativity of each subdiagram involved in it.
\begin{itemize}
\item $1$ commutes obviously since $S_n(\Delta)(\kappa_n)=(\kappa_n,\kappa_n)$.
\item $2$ commutes since
\begin{align*}
\left[\bigoplus_{p+q=n}S_p(\tilde{\sigma}^{q})\otimes S_q(\sigma^{p}_0)\right]\circ
\mathcal{AW}(\kappa_n,\kappa_n)
&=\sum_{p+q=n}S_p(\tilde{\sigma}^{q})\otimes S_q(\sigma^{p}_0)
(\tilde{\delta}^{q}\otimes\delta^{p}_0)\\
&=\sum_{p+q=n}\kappa_p\otimes\kappa_q.
\end{align*}
\item $3$ commutes by naturality of $\mathcal{EZ}$.
\item $4$ commutes by compatibility of $\mathcal{EZ}$ and $\mathcal{AW}$~\cite[I.4.b)]{rational homotopy}.
\item $5$ commutes by naturality of $(\mathcal{EZ}\otimes \mathcal{EZ})\circ (id\otimes\tau\otimes id)$.
\item $6$ commutes since
$$(\vert proj_1\vert,\vert
proj_2\vert)\circ\pi_n=(\pi_n\times\pi_n)\circ (id\times\tau\times
id)\circ (id\times\Delta).$$
\item $7$ commutes by naturality of $\mathcal{AW}$.
\item $8$ does not commute. But the two different maps coming from $8$ coincide on the
image of $(\mathcal{EZ}\otimes \mathcal{EZ})\circ (id\otimes\tau\otimes id)\circ [id\otimes \mathcal{AW}(\kappa_n,\kappa_n)]$.
Indeed this image is embedded in the image of
$$\sum_{p+q=n} S_*(id_{X_n}\times\tilde{\delta}^{q})\otimes S_*(id_{Y_n}\times\delta^{p}_0).$$
Now
$$\pi_p\circ(\tilde{d}^{q}\times\tilde{\sigma}^{q})\circ (id_{X_n}\times\tilde{\delta}^{q})
=\pi_p\circ(\tilde{d}^{q}\times id_{\Delta^{p}})
=\pi_n\circ(id_{X_n}\times\tilde{\delta}^{q}).$$
We have a similar formula for $Y_n$.
\end{itemize}
Finally, we have
$$(f\otimes f)\circ\left[\sum_{p,q} S_*(\tilde{d}^{q})\otimes S_*(d^{p}_0)\right]\circ\vert\mathcal{AW}\vert
=\mathcal{AW}\circ S_*((\vert proj_1\vert,\vert proj_2\vert)\circ f.$$\QED
\begin{lem}\cite[7.3.15]{Loday}\label{cyclique et lacets libres}
Let $X$ be a path connected pointed space. Then there is a natural homotopy equivalence
$\vert\Gamma\Omega X\vert\buildrel{\simeq}\over\rightarrow X^{S^{1}}$.
\end{lem}
\noindent{\bf Proof of Theorem \ref{Theoreme de Goodwillie}}
Applying Lemma \ref{echange cyclic et chaines singulieres} to the Moore loop space
$\Omega X$, Lemma \ref{realisation} to $\Gamma\Omega X$ and Lemma
\ref{cyclique et lacets libres} yield to the sequence of DGC quasi-isomorphisms:
$$\C S_*(\Omega X)=\vert\Gamma S_*(\Omega
X)\vert\buildrel{\simeq}\over\rightarrow \vert S_*(\Gamma\Omega X)\vert
\buildrel{\simeq}\over\rightarrow S_*(\vert \Gamma\Omega X\vert)
\buildrel{\simeq}\over\rightarrow S_*(X^{S^{1}}).$$\QED

\section{HAH models}
In order to compute the algebra structure of $HH^{*}S_*(\Omega X)$,
it is necessary to replace $S_*(\Omega X)$ by a smaller Hopf algebra.
Let's first remark that the cyclic bar construction preserves quasi-isomorphisms.
\begin{propriete}\cite[5.3.5]{Loday}(Compare \cite[4.3(iii)]{dga in topology})
\label{cyclic bar construction preserve quasi-isomorphisms}
Let $f:A\rightarrow B$ be a quasi-isomorphism of augmented DGA's.
If $\overline{A}$ and $\overline{B}$ are ${\Bbbk}$-semifree
then $\C(f):\C(A)\buildrel{\simeq}\over\rightarrow \C(B)$ is a quasi-isomorphism
of complexes.
\end{propriete}
Let $f$, $g$ : $A\rightarrow B$ be two morphisms of augmented DGA's.
A {\it derivation homotopy} from $f$ to $g$ is a morphism of graded modules of degree $+1$,
$h:A\rightarrow\overline{B}$
such that $d\circ h+h\circ d=f-g$
and $h(xy)=h(x)g(y)+(-1)^{\vert x\vert}f(x)h(y)$ for $x,y\in A$.
A derivation homotopy from $f$ to $g$ is denoted by $h:f\thickapprox g$.
We say that $f$ and $g$ are {\it homotopic} if there is a derivation homotopy between them.
\begin{lem}\label{cyclic bar construction et homotopie}
Let $f$, $g$ : $A\rightarrow B$ be two morphisms of augmented DGA's.
If $f$ and $g$ are homotopic then the morphisms of complexes
$\C(f)$, $\C(g)$ : $\C(A)\rightarrow \C(B)$ are chain homotopic.
\end{lem}
\pf Let $h$ be a derivation homotopy from $f$ to $g$.
By induction on the wordlength of the cyclic bar construction,
construct an explicit chain homotopy between $\C(f)$ and  $\C(g)$.\QED

A {\it Hopf algebra up to homotopy}, or HAH, is a DGA $K$ equipped with
two morphisms of DGA's $\Delta :K\rightarrow K\otimes K$ and $\varepsilon:K\rightarrow {\Bbbk}$
such that
$(\varepsilon\otimes id_K)\circ\Delta =id_K=(id_K\otimes \varepsilon)\circ\Delta$ (counitary exactly),
$(\Delta\otimes 1)\circ\Delta \thickapprox (1\otimes \Delta)\circ\Delta$ (coassociative up to homotopy) and
$\tau\circ\Delta\thickapprox\Delta$ (cocommutative up to homotopy).

Let $K$, $K'$ be two HAH's.
A morphism of augmented DGA's $f:K\rightarrow K'$ is a {\it HAH morphism} if
$\Delta f\thickapprox (f\otimes f)\Delta$ ($f$ commutes with the diagonals up to homotopy).
\begin{defin}\label{definition HAH model}
Let $X$ be a pointed topological space.
A {\it HAH model} for $X$ is a free chain algebra $(\TA V,\partial)$ equipped with a structure
of Hopf algebras up to homotopy and with a HAH quasi-isomorphism
$\Theta:(\TA V,\partial)\buildrel{\simeq}\over\rightarrow S_*(\Omega X)$.
\end{defin}
The existence of HAH models for any space is guarantied by the following two properties.
\begin{propriete}\cite[3.1]{dga in topology}\label{existence free model}
Let $A$ be a chain algebra. There exists a free chain algebra $(\TA V,\partial)$
and a quasi-isomorphism of augmented chain algebras
$$\Theta:(\TA V,\partial)\buildrel{\simeq}\over\longrightarrow A.$$
\end{propriete}
\begin{propriete}\cite[I.7 and II.1.11=II.2.11a)]{algebraic homotopy}
or \cite[3.6]{dga in topology}\label{lifting lemma}
Consider a quasi-isomorphism of augmented chain algebras $p:A\buildrel{\simeq}\over\longrightarrow B$
and a morphism of augmented chain algebras $g$ from a free chain algebra $(\TA V,\partial)$
to $B$:

$$\xymatrix{
&A\ar[d]^{p}_{\simeq}\\
(\TA V,\partial)\ar@{-->}[ur]^{f}\ar[r]_g
&B
}$$
Then
\begin{itemize}
\item[i)] there is a morphism of augmented chain algebras $f:(\TA V,\partial)\rightarrow A$
such that $p\circ f\thickapprox g$,
\item[ii)] moreover, any two such morphisms $f$ are homotopic.
\end{itemize}
\end{propriete}
Indeed by Property \ref{existence free model}, we obtain a quasi-isomorphism of augmented DGA's
$$\Theta:(\TA V,\partial)\buildrel{\simeq}\over\longrightarrow S_*(\Omega X).$$
Since $(\TA V,\partial)$ and $S_*(\Omega X)$ are ${\Bbbk}$-semifree
\footnote{We do not assume that ${\Bbbk}$ is a principal ideal domain~\cite[\S 2]{dga in topology}.},
$\Theta\otimes\Theta$ is a quasi-isomorphism.
By Property \ref{lifting lemma} i), we obtain a diagonal $\Delta_{\TA V}$ for $(\TA V,\partial)$
such that the following diagram of augmented DGA's commutes up to homotopy:
$$
\xymatrix{
(\TA V,\partial)\ar@{-->}[d]_{\Delta_{\TA V}}\ar[r]^{\simeq}_{\Theta}
&S_*(\Omega X)\ar[d]^{\Delta_{S_*(\Omega X)}}\\
(\TA V,\partial)\otimes (\TA V,\partial)\ar[r]^{\simeq}_{\Theta\otimes\Theta}
&S_*(\Omega X)\otimes S_*(\Omega X)
}$$
Since $S_*(\Omega X)$ is exactly a DGH and is cocommutative up to homotopy,
by Property \ref{lifting lemma} ii), $\Delta$ is counitary, coassociative and cocommutative up to
homotopy.
The diagonal $\Delta$ can be chosen to be strictly counitary~\cite[Lemma 5.4]{anick}.
\begin{theor}\label{cohomologie lacets libres et HAH}
Let $X$ be a path connected pointed space.
Let $(\TA V,\partial)$ be a HAH model for $X$.
There is a isomorphism of graded algebras $$HH^{*}(\TA V,\partial)\cong H^{*}(X^{S^{1}}).$$
\end{theor}
\pf Let $\Theta:(\TA V,\partial)\buildrel{\simeq}\over\rightarrow S_*(\Omega X)$ denote the HAH
quasi-isomorphism.
By Property \ref{cyclic bar construction preserve quasi-isomorphisms},
$\C(\Theta)$ is a quasi-isomorphism.
According Lemma \ref{cyclic bar construction et homotopie},
$$\C((\Theta\otimes\Theta)\circ\Delta_{\TA V})\thickapprox \C(\Delta_{S_*(\Omega X)}\circ\Theta).$$
So by composing with $AW$ and by applying Theorem \ref{Theoreme de Goodwillie},
the quasi-isomorphisms of chain complexes
$$\C(\TA V,\partial)\xrightarrow{\C(\Theta)} \C S_*(\Omega X)\buildrel{\simeq}\over\rightarrow S_*(X^{S^{1}})$$
commutes with the diagonals up to chain homotopy.\QED

\section{A smaller resolution than the bar resolution
}\label{smaller coalgebra}

The goal of this section is to replace the huge algebra up to homotopy
$\C(\TA V,\partial )^{\vee}$ by a smaller in order to be able to compute
the algebra $HH^{*}(\TA V,\partial)$.
When ${\Bbbk}$ is a field, Micheline Vigu\'e in \cite{Micheline} gives
a small complex $\left(({\Bbbk}\oplus sV)\otimes \T V,\delta\right)$
whose
homology is the vector space $HH_*(\TA V,\partial)$.
In fact, in this section, we show that over any commutative ring
${\Bbbk}$, this complex $\left(({\Bbbk}\oplus sV)\otimes
  \T V,\delta\right)$ is a strong deformation retract of $\C(\TA V,\partial)$.
\begin{defin}
Let $(Y,d)$ be a complex. A complex $(X,\partial)$ is a {\it strong deformation retract} of
$(Y,d)$ if there exist two morphisms of complexes $\nabla:(X,\partial)\hookrightarrow (Y,d)$,
$f:(Y,d)\twoheadrightarrow (X,\partial)$ and a chain homotopy $\Phi:(Y,d)\rightarrow (Y,d)$
such that
$f\nabla=id_X$ and $\nabla f-id_Y=d\Phi+\Phi d$.
The map $f$ is called the {\it projection} and the map $\nabla$ is called the {\it inclusion}. 
\end{defin}

We first consider the case where the differential $\partial$ on $(\TA V,\partial)$ is just obtained
by tensorization of the differential of a complex $V$ and is so
therefore homogeneous by wordlength.

Consider the tensor algebra $\TA V$ on a complex $V$.
Define the augmentation on $\TA V$ such that the augmentation ideal $\overline{\TA V}$ is
$$\T^{+}V=\oplus_{i\geq 1} V^{\otimes i}.$$
The bar resolution $\B(\TA V;\TA V;\TA V)$ contains a subcomplex
$(\T V\otimes({\Bbbk}\oplus sV)\otimes \T V,d_1+d_2)$, since
$$d_2(a\otimes sv\otimes b)=(-1)^{\vert a\vert}(av\otimes b-a\otimes vb).$$
\begin{proposition}\cite[Proposition 3.1.2]{Loday}\label{cas differentielle lineaire}
The $(\TA V,\TA V)$-bimodule $(\T V\otimes({\Bbbk}\oplus sV)\otimes \T V,d_1+d_2)$
is a strong deformation retract of the bar resolution $\B(\TA V;\TA V;\TA V)$.
\end{proposition}
\pf Define the projection $f:\B(\TA V;\TA V;\TA V)\twoheadrightarrow
\T V\otimes({\Bbbk}\oplus sV)\otimes \T V$ on its components
$f_n:\T V\otimes (s\T^{+}V)^{\otimes n}\otimes \T V\rightarrow \T V\otimes ({\Bbbk}\oplus sV)\otimes \T V$:

The map $f_0:\T V\otimes \T V\rightarrow \T V\otimes \T V$ is the identity map.

We define $f_1:\T V\otimes s\T^{+}V\otimes \T V\rightarrow \T V\otimes sV\otimes \T V$
by $$f_1(a[sv_1\cdots v_n]b)=\sum_{i=1}^{n}(-1)^{\vert v_1\cdots v_{i-1}\vert}
av_1\cdots v_{i-1}\otimes sv_i\otimes v_{i+1}\cdots v_nb$$
for $a,b\in \T V$, $v_1,\cdots,v_n\in V$ and $n\in \mathbb{N}^{*}$.

For $n\geq 2$, $f_n$ is the zero map.
An easy calculation shows that $d_2f_1=f_0d_2$.
Since $f_1(a[sa_1a_2]b)=f_1(a[sa_1]a_2b)+(-1)^{\vert a_1\vert} f_1(aa_1[sa_2]b)$,
$f_1d_2=0$. Therefore $f$ commutes with $d_2$ and is a morphism of complexes.

Of course, $f\nabla=id_{\T V\otimes({\Bbbk}\oplus sV)\otimes \T V}$.
The components $$\Phi_n:\T V\otimes (s\T^{+}V)^{\otimes n}\otimes \T V\rightarrow
\T V\otimes (s\T^{+}V)^{\otimes {n+1}}\otimes \T V$$ of the chain homotopy $\Phi$
are defined by:
\begin{eqnarray*}
\Phi_0&=&0,\\
\Phi_n(a[sa_1\vert\cdots\vert sa_{n-1}\vert sv]b)&=& 0,\\
\Phi_n(a[sa_1\vert\cdots\vert sa_{n-1}\vert sa_nv]b)& =
& -(-1)^{\varepsilon_n} a[sa_1\vert\cdots\vert sa_n\vert sv]b\\
&&+\Phi_n(a[sa_1\vert\cdots\vert sa_n]vb)
\end{eqnarray*}
for $a,b\in \T V$, $v\in V$ and $a_1,\cdots,a_n\in \T^+V$.
Recall that $\varepsilon_n=\vert a\vert+\vert sa_1\vert+\cdots+\vert sa_n\vert$.

By a double induction first on $n$ and then on the wordlength, check that
$d\Phi_n+\Phi_{n-1}d=\nabla f_n-id$, $n\in\mathbb{N}$.
At the beginning for $n=1$, use the formula
$f_1(a[sa_1v]b)=f_1(a[sa_1]vb)+(-1)^{\vert a_1\vert}aa_1\otimes sv\otimes b$.\QED

Consider now an augmented DGA $(\TA V,\partial)$ such that $\overline{\TA V}=\T^+V$.
The differential $\partial$ decomposes uniquely as a sum
$d_1+d_2+\cdots+d_i+\cdots$ of derivations satisfying
$d_i(V)\subset \T^iV=V^{\otimes i}$.
The differential $d_1$ is called the {\it linear part} of $d$.

To pass from the case $\partial=d_1$ to the general case, we'll use the
well-known perturbation Lemma.
For an abundant and recent bibliography, see \cite{Lambe} or \cite{Hess}.
\begin{theor}(Perturbation Lemma)
Let $(X,\partial)\build\leftrightarrows_\nabla^f (Y,d)\circlearrowleft\Phi$ be a strong deformation
retract of chain complexes satisfying $f\Phi=0$, $\Phi\nabla=0$ and $\Phi^2=0$.
Suppose moreover that this strong deformation retract is filtered:
there exist on $X$ and on $Y$ increasing filtrations bounded below
preserved by $\partial$, $d$, $f$, $\nabla$ and $\Phi$.
Consider a filtration lowering linear map $t:Y\rightarrow Y$ of degree $-1$
such that $d+t$ is a new differential on $Y$ (Such $t$ is called a
{\it perturbation}).
Then
\begin{align*}
\partial_\infty&=\partial+\sum_{k>0} f(t\Phi)^{k-1}t\nabla,\\
\nabla_\infty&=\nabla+\sum_{k>0} (\Phi t)^k\nabla,\\
f_\infty&=f+\sum_{k>0}f(t\Phi)^k,\\
\Phi_\infty&=\Phi+\sum_{k>0}(\Phi t)^k\Phi
\end{align*}
are well defined linear maps
and $(X,\partial_\infty)\build\leftrightarrows_{\nabla_\infty}^{f_\infty} (Y,d+t)
\circlearrowleft\Phi_\infty$ is a strong deformation
retract.
\end{theor}
By applying the Perturbation Lemma to Proposition \ref{cas differentielle lineaire}
we rediscover
\begin{theor}\cite[Th\'eor\`eme 1.4]{Micheline}\label{retract bar resolution}
Let $(\TA V,d)$ be a chain algebra.
Suppose that $V$ is a graded module concentrated in degree greater or equal than one.
Define the linear map of degree $+1$
\begin{eqnarray*}
S:\T V\otimes \T V&\rightarrow&\T V\otimes sV\otimes \T V\\
a\otimes v_1\cdots v_n&\mapsto&\sum_{i=1}^n(-1)^{\vert av_1\cdots v_{i-1}\vert}
av_1\cdots v_{i-1}\otimes sv_i\otimes
v_{i+1}\cdots v_n
\end{eqnarray*}
Consider the chain complex $(\T V\otimes ({\Bbbk}\oplus sV)\otimes \T V,D)$
where $$D\vert_{\T V\otimes \T V}=d,\qquad
D\vert_{\T V\otimes sV\otimes \T V}=\tilde{d_1}+d_2,$$
$$\tilde{d_1}(a\otimes sv\otimes b)=da\otimes sv\otimes b-S(a\otimes dv).b
-(-1)^{\vert av\vert}a\otimes sv\otimes db$$
$$\mbox{and}\quad d_2(a\otimes sv\otimes b)=(-1)^{\vert a\vert}(av\otimes b-a\otimes vb)
\mbox{ for }a,b\in \T V,\, v\in V.$$
Then $\left(\T V\otimes ({\Bbbk}\oplus sV)\otimes \T V,D\right)$ is a strong deformation retract
of the $(\TA V,\TA V)$-bimodule $\B(\TA V,\TA V,\TA V)$.
\end{theor}
\pf By Proposition \ref{cas differentielle lineaire}
$\left(\T V\otimes ({\Bbbk}\oplus sV)\otimes \T V,d_1+d_2\right)$
is a strong deformation retract of $\B\left((\TA V,d_1),(\TA V,d_1),(\TA V,d_1)\right)$
where $d_1$ denotes the linear part of $d$.
The anhilation conditions are satisfied:

$\Phi_1(a[sv]b)=0$ and $\Phi_0=0$, therefore $\Phi\nabla=0$.

The projection $f$ is null on $\T V\otimes (s\T^+V)^{\otimes\geq 2}\otimes \T V$
and $\Phi_0=0$. Therefore $f\Phi_n=0$ for $n\in\mathbb{N}$.

Since $\Phi_{n+1}\Phi_n(a[sa_1\vert\cdots\vert sa_nv]b)=
\Phi_{n+1}\Phi_n(a[sa_1\vert\cdots\vert sa_n]vb)$, by induction on wordlength
$\Phi_{n+1}\Phi_n=0$ for $n\geq 1$.

Let $n\in\mathbb{Z}$. An element $a[sa_1\vert\cdots\vert sa_n]b$ is
said to have a filtration degree $-n$ if and only if the sum of the wordlengths
of $a$, $a_1$, \dots, $a_n$ and $b$ is greater or equal than $n$.
The filtrations are bounded below since $V=V_{\geq 1}$.
The maps $\nabla$, $f$, $\Phi$, $d_1$ and $d_2$ respect wordlengths. Therefore
the strong deformation retract is filtered.
Define the perturbation $t$ to be equal to the differential of
$\B\left((\TA V,d),(\TA V,d),(\TA V,d)\right)$ minus the differential of
$\B\left((\TA V,d_1),(\TA V,d_1),(\TA V,d_1)\right)$.
Since $d_{\geq 2}=d-d_1$ increases wordlength by $1$ at least, $t$ is filtration lowering.

So finally we can apply the Perturbation Lemma and
$\left(\T V\otimes ({\Bbbk}\oplus sV)\otimes \T V,\partial_\infty\right)$
is a strong deformation retract of $\B\left((\TA V,d),(\TA V,d),(\TA V,d)\right)$.

The composite $t\Phi_n$ maps $\T V\otimes (s\T^+V)^{\otimes n}\otimes \T V$
into $\T V\otimes (s\T^+V)^{\otimes n+1}\otimes \T V$, $\Phi_0$ is null.
Therefore $f(t\Phi)^k=0$ for $k\geq 1$.
So $f=f_{\infty}$ (The projection is unchanged) and
$D:=\partial_{\infty}=\partial+ft\nabla=d_2+fd_1\nabla$ where
$d_1$ is the linear part of the the differential of
$\B\left((\TA V,d),(\TA V,d),(\TA V,d)\right)$.
Set $\tilde{d_1}:=fd_1\nabla$.
$$\tilde{d_1}(a\otimes sv\otimes b)=da\otimes sv\otimes b-(-1)^{\vert a\vert}
f_1(a\otimes sdv\otimes b)-(-1)^{\vert av\vert}a\otimes sv\otimes db.$$
$$S(a\otimes dv).b=(-1)^{\vert a\vert}f_1(a\otimes sdv\otimes b).$$\QED
\begin{cor}\cite[Th\'eor\`eme 1.5]{Micheline}\label{Homologie de Hochschild
d'une algebre libre}
Let $(\TA V,d)$ be a chain algebra such that $V=V_{\geq 1}$.
Define the linear map of degree $+1$
\begin{eqnarray*}
\overline{S}:\T V\otimes \T V & \rightarrow & sV\otimes \T V\\
v_1\ldots v_n\otimes a &\mapsto & \sum_{i=1}^n
(-1)^{\vert v_1\ldots v_{i-1}\vert\vert v_i\ldots v_n a\vert}
sv_i\otimes v_{i+1}\ldots v_nav_1\ldots v_{i-1}
\end{eqnarray*}
Consider the complex $\left(({\Bbbk}\oplus sV)\otimes \T V,\delta\right)$ where
$$\delta\vert_{\T V}=d,$$
$$\delta(sv\otimes a)=(-1)^{\vert a\vert\vert v\vert} 1\otimes av-1\otimes va
+(-1)^{\vert sv\vert}sv\otimes da-\overline{S}(dv\otimes a).$$
Then $\left(({\Bbbk}\oplus sV)\otimes \T V,\delta\right)$ is a strong deformation
retract of $\C(\TA V,d)$.
\end{cor}
\pf The linear maps
\begin{eqnarray*}
({\Bbbk}\oplus sV)\otimes \T V &\rightarrow&
\left(\T V\otimes ({\Bbbk}\oplus sV)\otimes \T V\right)
\otimes_{\TA V\otimes \TA V^{op}} \T V\\
\overline{v}\otimes b&\mapsto& 1\otimes\overline{v}\otimes 1\otimes b
\end{eqnarray*}
\begin{eqnarray*}
\left(\T V\otimes ({\Bbbk}\oplus sV)\otimes \T V\right)\otimes_{\TA V\otimes \TA V^{op}} \T V
&\rightarrow &({\Bbbk}\oplus sV)\otimes \T V\\
a\otimes\overline{v}\otimes a'\otimes b &\mapsto&
(-1)^{\vert a\vert (\vert\overline{v}\vert+\vert a'\vert+\vert b\vert)}
\overline{v}\otimes a'ba
\end{eqnarray*}
are inverse. The strong deformation retract given by Theorem~\ref{retract bar resolution}
is compatible with the structure of $(\TA V,\TA V)$-bimodule on $\B(\TA V;\TA V;\TA V)$.
To obtain the Corollary, we should tensor it by $\T V\otimes_{\TA V\otimes \TA V^{op}} -$
and then permute $\T V$ and $({\Bbbk}\oplus sV)$.
But it is equivalent and shorter to tensor by
$-\otimes_{\TA V\otimes \TA V^{op}} \T V$ and use the previous isomorphisms.\QED

Suppose that $(\TA V,\partial)$ is a HAH model of a path connected space $X$.
Using the inclusion $\nabla_{\infty}$ and the projection $f=f_{\infty}$
of the strong deformation retract
given by Corollary~\ref{Homologie de Hochschild d'une algebre libre},
it is now easy to transport the diagonal of $\C(\TA V,\partial)$,
denoted $\Delta_{\C(\TA V,\partial)}$, to
$\left(({\Bbbk}\oplus sV)\otimes \T V,\delta\right)$.
Define the diagonal of $\left(({\Bbbk}\oplus sV)\otimes \T V,\delta\right)$
simply as the composite $(f\otimes f)\circ
\Delta_{\C(\TA V,\partial)}\circ \nabla_{\infty}$.
Now $\left(({\Bbbk}\oplus sV)\otimes \T V,\delta\right)^{\vee}$ is an
algebra up to homotopy whose homology is isomorphic to $HH^{*}(\TA V,\partial)$
as graded algebras.
This algebra up to homotopy is the smallest that computes in general
the cohomology algebra of the free loop space $H^{*}(X^{S^{1}})$.
But the formula for the diagonal of $\left(({\Bbbk}\oplus sV)\otimes
  \T V,\delta\right)$ is very complicated: it involves in particular
the formula of the inclusion $\nabla_{\infty}$ given by the
Perturbation Lemma.

We will now limit ourself to two important cases where the HAH
structure on $(\TA V,\partial)$ is simple:
\begin{itemize}
\item The differential $\partial$ is the sum $d_1+d_2$ of only its
  linear part $d_1$ and its quadratic part $d_2$. The elements of $V$
are primitive: $\TA V$ is a primitively generated Hopf algebra. This will
be the object of Section~\ref{homologie de Hochschild et Cobar}.
\item The differential $\partial$ is equal to its linear part $d_1$
(hypothesis of Proposition~\ref{cas differentielle lineaire}).
The reduced diagonal $\overline{\Delta}$ of $\TA V$ is such that
$\overline{\Delta}(V)\subset V\otimes V$. This will
be the object of Section~\ref{free loop space on a suspension}.
\end{itemize}
\section{The isomorphism between $HH^{*}(\Omega C)$ and
$HH_*(C^{\vee})$}\label{homologie de Hochschild et Cobar}
Let $C$ be a coaugmented DGC. Denote by $\overline{C}$ the kernel of the counit.
The {\it cobar construction on $C$}, denoted $\Omega C$, 
is the augmented DGA $\left(\TA(s^{-1}\overline{C}),d_1+d_2\right)$
where $d_1$ and $d_2$ are the unique derivations determined by
$$d_1s^{-1}c=-s^{-1}dc\mbox{ and}$$
$$d_2s^{-1}c=\sum_i (-1)^{\vert x_i\vert} s^{-1}x_i\otimes s^{-1}y_i,\; c\in\overline{C}$$
where the reduced diagonal $\displaystyle\overline{\Delta}c=\sum_i x_i\otimes y_i$. We follow the sign convention of \cite{Adams' cobar equivalence}.
\begin{theor}\cite[Theorem A]{Jones et McCleary}\cite[Theorem II]{Steve et Micheline}\label{formalite et homologie de Hochschild}
Consider a coaugmented DGC $C$ $\Bbbk$-free of finite type such that $C={\Bbbk}\oplus C_{\geq 2}$.
Then there is a natural isomorphism of graded modules
$$HH^{*}(\Omega C)\cong HH_*(C^{\vee}).$$
\end{theor}
We give again the proof of Jones and McCleary since we want to check
carefully the signs. We also need to explicit the isomorphism in order
to transport later the algebra structure.
We remark that already at the level of complexes, there is a
quasi-isomorphism from $\C(\Omega C)^{\vee}$ to $\C(C^{\vee})$.

Before giving the proof, we need to give the signs convention used in
this paper.
Let $f:V\rightarrow V'$ and $g:W\rightarrow W'$ be two linear maps
then $f\otimes g:V\otimes W\rightarrow V'\otimes W'$ is the linear map
given by
$$(f\otimes g)(v\otimes w)=(-1)^{\vert g\vert\vert v\vert}f(v)\otimes
g(w).$$
Therefore if $f':V'\rightarrow V"$ and $g':W'\rightarrow W"$ are two
other
linear maps then
$$(f'\otimes g')\circ(f\otimes g)=(-1)^{\vert f\vert\vert g'\vert}
(f'\circ g')\otimes (f\circ g).$$
Let $\varphi:M\rightarrow N$ be a linear map.
If $f\in\mbox{Hom}(N,{\Bbbk})$ then $$\varphi^{\vee}(f)=(-1)^{\vert
  f\vert\vert\varphi\vert}f\circ\varphi .$$
In particular, if $(M,d)$ is a complex, the dual complex is
$(M^{\vee},d^{\vee})$.
Let $\Psi:N\rightarrow Q$ be another linear map.
Then $(\Psi\circ\varphi)^{\vee}=(-1)^{\vert\varphi\vert\vert\Psi\vert}
\varphi^{\vee}\circ\Psi^{\vee}$.
\pf We apply Corollary \ref{Homologie de Hochschild d'une algebre libre}
when the chain algebra $(\TA V,d)$ is the cobar $\Omega C$. We obtain a strong deformation
retract of the cyclic bar construction $\C(\Omega C)$ of the form
$(C\otimes\Omega C,\delta)$. The differential $\delta$ is given by
\begin{eqnarray*}
\delta a&=&d_\Omega a,\; a\in\Omega C,\\
\delta (c\otimes a)&=&dc\otimes a+(-1)^{\vert c\vert} c\otimes da\\
&&-1\otimes (s^{-1}c)a-(-1)^{\vert x_i\vert} x_i\otimes (s^{-1}y_i)a\\
&&+(-1)^{\vert a\vert\vert s^{-1}c\vert}1\otimes a s^{-1}c
+(-1)^{(\vert a\vert +\vert y_i\vert)\vert s^{-1}x_i\vert}y_i\otimes a s^{-1}x_i.
\end{eqnarray*}
Therefore $(C\otimes\Omega C,\delta)$ is the complex $(C\otimes \T s^{-1}\overline{C},d_1+d_2)$ where $d_1$ is just obtained by tensorization and
$d_2:C\otimes (s^{-1}\overline{C})^{\otimes n-1}\rightarrow
C\otimes (s^{-1}\overline{C})^{\otimes n}$ is the sum of $n+1$ terms
$\delta_0$, $\delta_1$, \dots, $\delta_n$ given by
$$\delta_0=-\left[(C\otimes s^{-1})\circ\Delta\right]\otimes (s^{-1}\overline{C})^{\otimes n-1},$$
$$\delta_i=C\otimes (s^{-1}\overline{C})^{\otimes i-1}\otimes
\left[(s^{-1}\otimes s^{-1})\circ\Delta\circ s\right]\otimes
(s^{-1}\overline{C})^{\otimes n-1-i},\;1\leq i\leq n-1$$
$$\mbox{and }\delta_n=
\left[C\otimes (s^{-1}\overline{C})^{\otimes n-1}\otimes s^{-1}\right]\circ
\tau_{C,C\otimes (s^{-1}\overline{C})^{\otimes n-1}}\circ
\left[\Delta\otimes (s^{-1}\overline{C})^{\otimes n-1}\right].$$
Let $A$ denote the augmented DGA $C^{\vee}$.
The differential $d_2:A\otimes (s\overline{A})^{\otimes n}\rightarrow
A\otimes (s\overline{A})^{\otimes n-1}$ is also the sum of $n+1$ terms
$d_0$, $d_1$, \dots, $d_n$ (compare to Example~\ref{Bar construction cyclique simpliciale})
given by
$$d_0=\left[\mu\circ(A\otimes s^{-1})\right]\otimes (s\overline{A})^{\otimes n-1},$$
$$d_i=A\otimes (s\overline{A})^{\otimes i-1}\otimes
\left[s\circ\mu\circ(s^{-1}\otimes s^{-1})\right]\otimes (s\overline{A})^{n-i-1}
,\;1\leq i\leq n-1$$
$$\mbox{and }d_n=-\left[\mu\otimes (s\overline{A})^{\otimes n-1}\right]\circ
\tau_{A\otimes (s\overline{A})^{\otimes n-1},A}\circ
\left[A\otimes (s\overline{A})^{\otimes n-1}\otimes s^{-1}\right].$$
The isomorphism $\Theta:s(\overline{C}^{\vee})\buildrel{\cong}\over
\rightarrow (s^{-1}\overline{C})^{\vee}$ is such that $(s^{-1})^{\vee}\circ\Theta=s^{-1}$
\cite[p. 276]{universal enveloping algebras}.
For any two complexes $V$ and $W$, the map
$\Phi:V^{\vee}\otimes W^{\vee}\rightarrow (V\otimes W)^{\vee}$ given
by
$\Phi(f\otimes g)=\mu_{\Bbbk}\circ(f\otimes g)$
is a morphism of complexes and is associative, commutative, natural with respect
to linear maps of any degree.
Therefore the composite
$$A\otimes (s\overline{A})^{\otimes n}\buildrel{A\otimes\Theta^{\otimes n}}\over
\longrightarrow C^{\vee}\otimes \left[(s^{-1}\overline{C})^{\vee}\right]^{\otimes n}
\buildrel{\Phi}\over\rightarrow
\left[C\otimes(s^{-1}\overline{C})^{\otimes n}\right]^{\vee}$$
commutes with $d_i$ and $\delta_i^{\vee}$ for $0\leq i\leq n$.
So finally $$\Phi\circ\left[A\otimes \T(\Theta)\right]:\C(A)
\buildrel{\cong}\over\rightarrow (C\otimes\Omega C,\delta)^{\vee}$$ is
an isomorphism of complexes.\QED

Let $V$ be a graded module.
The tensor algebra $\TA V$ can be made into a cocommutative Hopf algebra by requiring
the elements of $V$ to be primitive \cite[0.5 (10)]{Tanre}.
We will called the resulting diagonal, the {\it shuffle diagonal}.
Dually the tensor coalgebra $\TC V$ equipped with the {\it shuffle product}
\label{shuffle product}
is a commutative Hopf algebra. The shuffle product is defined by
$$[v_1|\ldots|v_p]\cdot[v_{p+1}|\ldots|v_{p+q}]=
\sum_{\sigma} \sigma\cdot[v_1|\ldots|v_{p+q}]$$
where the sum is taken over the $(p,q)$-shuffles $\sigma$ and a permutation $\sigma$
acts on $[v_1|\ldots|v_{p+q}]$ by permuting the factors with appropriate signs
\cite[Appendix]{universal enveloping algebras} or \cite[0.5 (8)]{Tanre}.
Suppose that $V$ is ${\Bbbk}$-free of finite type and $V=V_{\geq 1}$.
Then the map $\Phi:\TC(V^{\vee})\buildrel{\cong}\over\rightarrow (\TA V)^{\vee}$
is an isomorphism of Hopf algebras.

Let $C$ be a cocommutative coaugmented DGC. Then the cobar
$\Omega C=(\TA(s^{-1}\overline{C}),d_1+d_2)$ equipped with the shuffle diagonal
is a DGH \cite[0.6 (2)]{Tanre}.
Dually, let $A$ be an augmented CDGA.
Consider the multiplication on $A\otimes \TC(s\overline{A})\otimes A$
obtained by tensorizing the multiplication of $A$ and the shuffle product
of $\TC(s\overline{A})$.
Then the bar resolution of $A$, $\B(A;A;A)=(A\otimes \TC(s\overline{A})\otimes A,d_1+d_2)$
is a CDGA. The cyclic bar construction $\C(A)=(A\otimes \TC(s\overline{A}),d_1+d_2)$
is also a CDGA.
Therefore the Hochschild homology of a CDGA $A$, $HH_*(A)$,
has a natural structure of commutative graded algebra \cite[4.2.7]{Loday}.
The reduced bar construction of $A$, $\B(A)=(\TC(s\overline{A}),d_1+d_2)$
is a commutative DGH \cite[0.6 (1)]{Tanre}.
\begin{theor}\label{shc formalite et homologie de Hochschild}
Under the hypothesis of Theorem \ref{formalite et homologie de Hochschild},
if $C$ is cocommutative
then the isomorphism $HH^{*}(\Omega C)\cong HH_*(C^{\vee})$
is an isomorphism of commutative graded algebras.
\end{theor}
\begin{propriete}\label{diagonale de la bar sur primitifs}
Let $K$ be a graded Hopf algebra.
Consider $\mbox{Ker}\overline{\Delta}$ the primitive elements of $K$
and the graded coalgebra $K\otimes \TC(s\mbox{Ker}\overline{\Delta})\otimes K$
obtained by tensorization. Then the canonical map
$K\otimes \TC(s\mbox{Ker}\overline{\Delta})\otimes K\rightarrow \B(K;K;K)$
is a morphism of graded coalgebras.
\end{propriete}
\noindent{\bf Proof of Theorem \ref{shc formalite et homologie de
    Hochschild}}

When restricted to conilpotent coaugmented DGC's,
the cobar construction $\Omega$ is a left adjoint functor to the bar construction $\B$
\cite[Proposition 2.11]{Adams' cobar equivalence}.
By Formula~\ref{formule du relevement coalgebre colibre}, the adjunction map
$\sigma_C:C\buildrel{\simeq}\over\rightarrow \B\Omega C$ is
given by
$$
\sigma_C(c)=\sum_{i=0}^{+\infty}\sum [ss^{-1}c_1\vert\cdots\vert ss^{-1}c_{i+1}],
\;c\in\overline{C}
$$
where the iterated reduced diagonal $\overline{\Delta}^{i}c=\sum
c_1\otimes\dots\otimes c_{i+1}$.
We consider now the inclusion map
$\nabla_{\infty}$ of the strong deformation retract given by Theorem
\ref{retract bar resolution} when the chain algebra $(\TA V,d)$ is the cobar $\Omega C$.
A simple computation shows that
$\nabla_{\infty}$ is $\Omega C\otimes\sigma_C\otimes\Omega C$, the tensor
product of the identity maps and the adjunction map.
Now $\sigma_C$ is a morphism of coalgebras,
$\mbox{Im}\sigma_C\subset \TC(ss^{-1}\overline{C})$ and
$s^{-1}\overline{C}\subset\mbox{Ker}\overline{\Delta}$.
Therefore by Property \ref{diagonale de la bar sur primitifs},
the coalgebra $\Omega C\otimes C\otimes\Omega C$ obtained by tensorization
is a sub DGC of $\B(\Omega C;\Omega C;\Omega C)$.
After tensorizing by $-\otimes_{\Omega C\otimes\Omega C^{op}}\Omega C$ and
dualizing, we obtain the natural DGA quasi-isomorphism
$$\C(C^{\vee})\buildrel{\cong}\over\rightarrow
(C\otimes\Omega C,\delta)^{\vee}\buildrel{\simeq}\over\leftarrow
\C(\Omega C)^{\vee}.$$\QED

\section{The free loop space on a suspension}\label{free loop space on
  a suspension}
In this section, we show how to compute the cohomology algebra
of the free loop space on any suspension,
$H^{*}\left((\Sigma X)^{S^{1}}\right)$.
And even better, from the DGA $S^{*}(X)$ or any DGA weakly equivalent,
we construct an DGA weakly equivalent to the DGA
$S^{*}\left((\Sigma X)^{S^{1}}\right)$.

First we introduce some terminology.
Let $C$ be a coaugmented DGC. The composite
$C\buildrel{\Delta_C}\over\rightarrow C\otimes C\hookrightarrow
\TA\overline{C}\otimes \TA\overline{C}$ extends to an unique morphism
of augmented DGA's
$$\Delta_{\TA\overline{C}}:\TA\overline{C}\rightarrow \TA\overline{C}\otimes
\TA\overline{C}.$$
This DGH structure on the tensor algebra $\TA\overline{C}$ is called
the {\it Hopf algebra structure obtained by tensorization of the
coalgebra $C$}.
Dually, let $A$ be an augmented DGA.
The composite $\TC\overline{A}\otimes \TC\overline{A}\twoheadrightarrow
A\otimes A\buildrel{\mu_A}\over\rightarrow A$ lifts to an unique morphism
of coaugmented DGC's
$$\mu_{\TC\overline{A}}:\TC\overline{A}\otimes \TC\overline{A}\rightarrow
\TC\overline{A}.$$
This DGH structure on the tensor coalgebra $\TC\overline{A}$ is called
the {\it Hopf algebra structure obtained by tensorization of the
algebra $A$}.
Using formula~\ref{formule du relevement coalgebre colibre},
we see that the product $\mu_{\TC\overline{A}}$ of two elements
$[a_1|\cdots|a_p]$ and $[b_1|\cdots|b_q]$ admits the following
description:

A sequence
$\sigma=\left( (0,0)=(x_0,y_0),(x_1,y_1),\dots,(x_n,y_n)=(p,q) \right)$
defined by $$(x_i,y_i)=
\begin{cases}
(x_{i-1}+1,y_{i-1}) &
\text{or}\\
(x_{i-1},y_{i-1}+1) &
\text{or}\\
(x_{i-1}+1,y_{i-1}+1),  
\end{cases}
$$
is called {\it a step by step path} of length $n$ from $(0,0)$ to $(p,q)$.
To any step by step path $\sigma$ of length $n$,
we associate $c_\sigma=[c_1|\cdots|c_n]\in\overline{A}^{\otimes n}$
by the rule
\[c_i=
\begin{cases}
a_{x_i} &
\text{if $(x_i,y_i)=(x_{i-1}+1,y_{i-1})$, $i^{th}$ step is toward right,}\\
b_{y_i} &
\text{if $(x_i,y_i)=(x_{i-1},y_{i-1}+1)$, $i^{th}$ step is toward up,}\\
\mu_A(a_{x_i}\otimes b_{y_i}) &
\text{if $(x_i,y_i)=(x_{i-1}+1,y_{i-1}+1)$, $i^{th}$ step in diagonal.}  
\end{cases}
\]
Then a straightforward computation establishes
\addtocounter{theor}{1}
\begin{equation}\label{produit cohomologie lacets sur suspension}
\mu_{\TC\overline{A}}([a_1|\cdots|a_p]\otimes
  [b_1|\cdots|b_q])=\sum_\sigma \pm c_\sigma
\end{equation} 
where the sum is taken over all the step by step paths $\sigma$ from $(0,0)$
to $(p,q)$ and where $\pm$ is the
sign obtained with Koszul rule by mixing the
$a_1,\cdots,a_p$ and the $b_1,\cdots,b_q$.
In particular, when the product of $A$ is trivial, the product
$\mu_{\TC\overline{A}}$ is the shuffle product
considered page~\pageref{shuffle product}.

Suppose that the DGC $C$ is ${\Bbbk}$-free of finite type and such that
$C={\Bbbk}\oplus C_{\geq 1}$.
Then the map
$$\Phi:\TC(\overline{C^{\vee}})\buildrel{\cong}\over\rightarrow
\left(\TA\overline{C}\right)^{\vee}$$ is a DGH isomorphism.

The starting observation of this section is the following consequence of
Bott-Samelson Theorem (see \cite[7.1]{article} for details).
\begin{lem}\label{Theoreme de Bott-Samelson}
Let $X$ be a path connected space.
\begin{itemize}
\item[i)] Consider the Hopf algebra structure on
$\TA\overline{S_*(X)}$ obtained by tensorization of the coalgebra
$S_*(X)$. Then there is a natural DGH quasi-isomorphism
$$\TA\overline{S_*(X)}\buildrel{\simeq}\over\rightarrow
S_*(\Omega\Sigma X).$$
\item[ii)] Suppose that
$H_*(X)$ is ${\Bbbk}$-free.
Consider the Hopf algebra structure on $\TA H_+(X)$ obtained by
tensorization on the coalgebra $H_*(X)$.
Then there is a HAH quasi-isomorphism
$$\Theta_X:\TA H_+(X)\buildrel{\simeq}\over\rightarrow
S_*(\Omega\Sigma X).$$
\end{itemize}
\end{lem}
We can choose $\Theta_X$ to be natural in homology and so natural after passing to
homotopy of algebras (Property \ref{lifting lemma} i)).
\begin{defin}\cite[2.2.11]{Parhizgar}
Let $R$ be a graded algebra. Let $M$ be a $(R,R)$-bimodule.
The graded module $R\oplus M$, product of $R$ and of $M$ equipped
with the multiplication
$$(r_1,m_1)(r_2,m_2)=(r_1r_2,r_1\cdot m_2+m_1\cdot r_2)$$
is a graded algebra called the {\it trivial extension of $R$ by $M$}.
\end{defin}
\begin{lem}\label{suspension partie algebrique}
Let $C$ be a chain coalgebra ${\Bbbk}$-free of finite type such that $C={\Bbbk}\oplus C_{\geq 1}$.
Denote by $A$ the cochain algebra dual of $C$.
Consider the Hopf algebra structures on $\TA\overline{C}$ and on
$\TC\overline{A}$ obtained by tensorization of the coalgebra $C$
and the algebra $A$.
Define a structure of $(\TC\overline{A},\TC\overline{A})$-bimodule
on $s^{-1}\overline{A}\otimes \T\overline{A}$ by
\begin{align*}
\begin{split}
(s^{-1}a\otimes m) \cdot a_1\dots a_n=&
-s^{-1}a\otimes \mu_{\TC\overline{A}}(m\otimes a_1\dots a_n)\\
&-(-1)^{\vert a_n\vert (\vert m\vert+\vert a_1\dots a_{n-1}\vert )}
s^{-1}\mu_A(a\otimes a_n)\otimes
\mu_{\TC\overline{A}}(m\otimes a_1\dots a_{n-1})\\
\intertext{and}
a_1\dots a_n \cdot (s^{-1}a\otimes m)=&
-(-1)^{\vert s^{-1} a\vert\vert a_1\dots a_n\vert}
s^{-1}a\otimes \mu_{\TC\overline{A}}(a_1\dots a_n\otimes m)\\
&-(-1)^{\vert a_1\vert+\vert s^{-1} a\vert\vert a_2\dots a_n\vert}s^{-1}\mu_A(a_1\otimes a)\otimes \mu_{\TC\overline{A}}(a_2\dots
a_n\otimes m)
\end{split}
\end{align*}
for $a\in\overline{A}$, $m\in \T\overline{A}$ and
$a_1\dots a_n\in \TC\overline{A}$.
Consider the cochain complex $({\Bbbk}\oplus s^{-1}\overline{A})$
equipped with the trivial product.
Then the cyclic bar construction on
$({\Bbbk}\oplus s^{-1}\overline{A})$,
$\C({\Bbbk}\oplus s^{-1}\overline{A})=\TC\overline{A}\oplus
(s^{-1}\overline{A}\otimes \T\overline{A})$, equipped with the product
of the trivial extension of $\TC\overline{A}$ by $s^{-1}\overline{A}\otimes \T\overline{A}$ is a cochain algebra equipped with a natural
DGA quasi-isomorphism
$$\C({\Bbbk}\oplus s^{-1}\overline{A})\buildrel{\simeq}\over
\twoheadleftarrow \C(\TA\overline{C})^{\vee}.$$
\end{lem}
\begin{propriete}\label{diagonal restreinte bar construction}
Let $K$ be a graded Hopf algebra.
The diagonal of the coalgebra $\B(K;K;K)$ restricted to
$K\otimes ({\Bbbk}\oplus s\overline{K})\otimes K$
is the $(K,K)$-linear map given by
$$\Delta [sx]=[sx]\otimes []+[sy]\otimes z[]+
(-1)^{\vert y\vert}[]y\otimes [sz]+[]\otimes [sx],
x\in\overline{K}$$
where the reduced diagonal
$\overline{\Delta} x=\sum y\otimes z$.
\end{propriete}
\noindent{\bf Proof of Lemma \ref{suspension partie algebrique}}
The tensor algebra $\TA\overline{C}$
is equal as DGA to the cobar on the DGC ${\Bbbk}\oplus s\overline{C}$
with trivial coproduct.
Therefore by Theorem~\ref{formalite et homologie de Hochschild},
we get immediatly a natural quasi-isomorphism of cochain complexes
$$\C({\Bbbk}\oplus s^{-1}\overline{A})\buildrel{\simeq}\over\leftarrow
\C(\TA\overline{C})^{\vee}.$$
But we have to remember how this morphism decomposes in order to
transport the multiplication from $\C(\TA\overline{C})^{\vee}$ to
$\C({\Bbbk}\oplus s^{-1}\overline{A})$.

Since the differential on $\TA\overline{C}$ is only linear,
by Proposition~\ref{cas differentielle lineaire},
the canonical inclusion
$$\T\overline{C}\otimes ({\Bbbk}\oplus s\overline{C})\otimes \T\overline{C}
\buildrel{\simeq}\over\hookrightarrow
\B(\TA\overline{C};\TA\overline{C};\TA\overline{C})$$
is a quasi-isomorphism of complexes.
Since the reduced diagonal of $\TA\overline{C}$,
$\overline{\Delta}_{\TA\overline{C}}$
embeds $\overline{C}$ into $\overline{C}\otimes \overline{C}$,
by Property~\ref{diagonal restreinte bar construction},
$\T\overline{C}\otimes ({\Bbbk}\oplus s\overline{C})\otimes \T\overline{C}$
is a sub coalgebra of
$\B(\TA\overline{C};\TA\overline{C};\TA\overline{C})$. By tensorizing by
$-\otimes_{\TA\overline{C}\otimes \TA\overline{C}} \TA\overline{C}$,
we obtain the DGC $({\Bbbk}\oplus s\overline{C})\otimes \T\overline{C}$
with differential given by Corollary~\ref{Homologie de Hochschild
d'une algebre libre} and diagonal given by
\begin{align*}
\Delta(sx\otimes c)=&sx\otimes c'\otimes 1\otimes c''
+(-1)^{\vert c\vert\vert z\vert}sy\otimes c'\otimes 1\otimes c''z\\
&+(-1)^{\vert y\vert+\vert c'\vert\vert sz\vert}
1\otimes yc'\otimes sz\otimes c"
+(-1)^{\vert c'\vert\vert sx\vert}1\otimes c'\otimes sx\otimes c"
\end{align*}
for $x\in\overline{C}$, $c\in \T\overline{C}$ and where
the reduced diagonal $\overline{\Delta}x=\sum y\otimes z$
and the unreduced diagonal $\Delta c=\sum c'\otimes c"$.
The canonical inclusion
$$({\Bbbk}\oplus s\overline{C})\otimes \T\overline{C}
\buildrel{\simeq}\over\hookrightarrow \C(\TA\overline{C})$$
is a DGC quasi-isomorphism.

In order to dualize (for details, review the proof of
Theorem~\ref{formalite et homologie de Hochschild}), we see that the
diagonal
on $({\Bbbk}\oplus s\overline{C})\otimes \T\overline{C}$
is the sum of three terms,
$$\Delta_{\TA\overline{C}}:\T\overline{C}\rightarrow \T\overline{C}\otimes
\T\overline{C},$$
$$\Delta_1:s\overline{C}\otimes \T\overline{C}\rightarrow
s\overline{C}\otimes \T\overline{C}\otimes \T\overline{C}$$
and $$\Delta_2:s\overline{C}\otimes \T\overline{C}\rightarrow
\T\overline{C}\otimes s\overline{C}\otimes \T\overline{C}.$$
The first term $\Delta_{\TA\overline{C}}$ is just the diagonal
of $\TA\overline{C}$.
The second term $\Delta_1$ is the composite
$$(s\otimes \T\overline{C}\otimes\mu_{\TA\overline{C}})\circ
(C\otimes \T\overline{C}\otimes \T\overline{C}\otimes i)\circ
(C\otimes\tau_{C,\T\overline{C}\otimes \T\overline{C}})\circ
(\Delta_C\otimes\Delta_{\TA\overline{C}})\circ
(s^{-1}\otimes \T\overline{C})$$
where $i$ denotes the inclusion $C\hookrightarrow \T\overline{C}$.
The third term $\Delta_2$ is the composite
$$(\mu_{\TA\overline{C}}\otimes s\otimes \T\overline{C})\circ
(i\otimes\tau_{C,\T\overline{C}}\otimes \T\overline{C})\circ
(\Delta_C\otimes\Delta_{\TA\overline{C}})\circ
(s^{-1}\otimes \T\overline{C}).$$
Therefore the product on
$({\Bbbk}\oplus s^{-1}\overline{A})\otimes \T\overline{A}$
is the sum of three terms: the product $\mu_{\TC\overline{A}}$ of
$\TC\overline{A}$, the dual of $\Delta_1:
s^{-1}\overline{A}\otimes \T\overline{A}\otimes \T\overline{A}\rightarrow
s^{-1}\overline{A}\otimes \T\overline{A}$
 and the dual of $\Delta_2:
\T\overline{A}\otimes s^{-1}\overline{A}\otimes
\T\overline{A}\rightarrow
s^{-1}\overline{A}\otimes \T\overline{A}$.

Explicitly the product is given by
\begin{align*}
\begin{split}
(1\otimes m).(1\otimes m')=&1\otimes \mu_{\TC\overline{A}}(m\otimes m'),\\
(s^{-1}a\otimes m).(1\otimes a_1\dots a_n)=&
-s^{-1}a\otimes \mu_{\TC\overline{A}}(m\otimes a_1\dots a_n)\\
&-\pm s^{-1}\mu_A(a\otimes a_n)\otimes
\mu_{\TC\overline{A}}(m\otimes a_1\dots a_{n-1}),\\
(1\otimes a_1\dots a_n).(s^{-1}a\otimes m)=&
-\pm s^{-1}a\otimes \mu_{\TC\overline{A}}(a_1\dots a_n\otimes m)\\
&-\pm s^{-1}\mu_A(a_1\otimes a)\otimes \mu_{\TC\overline{A}}(a_2\dots
a_n\otimes m)\\
\intertext{and}
(s^{-1}a\otimes m).(s^{-1}a'\otimes m')=&0.
\end{split}
\end{align*}
for $a$, $a'\in\overline{A}$, $m$, $m'$ and $a_1\dots a_n\in
\T\overline{A}$
and where $\pm$ are the signs obtained exactly by the Koszul sign
convention.\QED
\begin{theor}\label{modele des lacets libres sur suspension}
Let $X$ be a path connected space.
If $S_*(X)$ is weakly equivalent as ${\Bbbk}$-free chain coalgebra to
a chain coalgebra $C$ ${\Bbbk}$-free of finite type such that
$C={\Bbbk}\oplus C_{\geq 1}$.
Then the singular cochains on the free loop spaces on the suspension
of $X$, $S^{*}((\Sigma X)^{S^{1}})$ is weakly equivalent
as cochain algebras to the cyclic bar construction
$\C({\Bbbk}\oplus s^{-1}\overline{C^{\vee}})$ equipped with the product
of the trivial extension given by Lemma \ref{suspension partie algebrique}.
\end{theor}
\begin{rem}\label{DGC modele}
If X is a finite simply connected CW-complex then $X$ satisfies the
hypothesis
of Theorem \ref{modele des lacets libres sur suspension}.
Indeed, the Adams-Hilton model of $X$ denoted ${\cal A}(X)$ is
a free chain algebra $(\TA V,\partial)$ equipped with a quasi-isomorphism
of augmented chain algebras
$(\TA V,\partial)\buildrel{\simeq}\over\rightarrow S_*(\Omega X)$ and 
such that the complex of indecomposables $(V,d_1)$ is the
desuspension of the reduced cellular chain complex of $X$.
Now the bar construction $\B S_*(\Omega X)$ is weakly equivalent as
${\Bbbk}$-free chain coalgebras to $S_*(X)$.
Therefore we can take $C=B{\cal A}(X)$.
\end{rem}
\begin{ex}
$S^{d}$, $d\geq 1$. If $d\geq 2$ using Remark \ref{DGC modele},
$S_*(S^{d})$ is weakly DGC equivalent to $\B\Omega H_*(S^{d})$,
therefore to $H_*(S^{d})$. By \cite[7.3]{dga in topology},
there is a DGH quasi-isomorphism
$H_*(S^{1})\buildrel{\simeq}\over\rightarrow S_*(S^{1})$.
So by Theorem \ref{modele des lacets libres sur suspension}, as DGA
$$S^{*}\left((S^{d+1})^{S^{1}}\right)\sim \C H^{*}(S^{d+1})
=E(s^{-1}v)\otimes \TC v,d_2$$
where $v$ is an element of degree $d$.
If $d$ is even 
in ${\Bbbk}$, as DGA
$$S^{*}\left((S^{d+1})^{S^{1}}\right)\sim E(s^{-1}v)\otimes \Gamma
v,0$$
and $$H^{*}\left((S^{d+1})^{S^{1}}\right)\cong H^{*}(S^{d+1})\otimes
H^{*}(\Omega S^{d+1})$$ as graded algebras. 
We suppose now that $d$ is odd. By dualization, $\TC v\cong Ev\otimes
\TC(v^{2})$ as graded algebras (James-Toda).
So as cochain algebras $$S^{*}\left((S^{d+1})^{S^{1}}\right)\sim
E(s^{-1}v)\otimes Ev\otimes \Gamma(v^{2}),
d_2\gamma_k(v^{2})=2(s^{-1}v)v\gamma_{k-1}(v^{2}), k\geq 1.$$
Therefore, over any commutative ring ${\Bbbk}$, the graded algebra
$H^{*}\left((\Sigma S^{d})^{S^{1}}\right)$ is the module
$${\Bbbk}\oplus
(_2{\Bbbk})\Gamma^{+}(v^{2})\oplus{\Bbbk}v.\Gamma(v^{2})
\oplus{\Bbbk}(s^{-1}v).\Gamma(v^{2})\oplus(\frac{{\Bbbk}}{2{\Bbbk}})v.(s^{-1}v).\Gamma(v^{2})$$
equipped with the obvious products.
In particular, if $\frac{1}{2}\in{\Bbbk}$, all the products are trivial.
\end{ex}
\begin{ex}\label{Comparaison espace projectif et wedge
de spheres} Comparaison of $\mathbb{CP}^{d}$ and $S^{2}\vee\dots\vee
  S^{2d}$, $d\geq 1$.
The Adams-Hilton model of $\mathbb{CP}^{d}$ is $\Omega
H_*(\mathbb{CP}^{d})$.
Therefore (Remark \ref{DGC modele}), $S_*(\mathbb{CP}^{d})$ is weakly
DGC
equivalent to $H_*(\mathbb{CP}^{d})$.
So by Theorem \ref{modele des lacets libres sur suspension},
as cochain algebras
$$S^{*}\left((\Sigma \mathbb{CP}^{d})^{S^{1}}\right)\sim
\C(H^{*}(\Sigma \mathbb{CP}^{d})).$$
Similarly as cochain algebras
$$S^{*}\left((S^{3}\vee\dots\vee S^{2d+1})^{S^{1}}\right)\sim
\C(H^{*}(S^{3}\vee\dots\vee S^{2d+1})).$$
The cochain algebras $\C(H^{*}(\Sigma \mathbb{CP}^{d}))$ and
$\C(H^{*}(S^{3}\vee\dots\vee S^{2d+1}))$ have the same underlying
cochain complex. But the non commutative product
on $\C(H^{*}(\Sigma \mathbb{CP}^{d}))$ given by Lemma
\ref{suspension partie algebrique} and 
using formula \ref{produit cohomologie lacets sur suspension}, is
far more complicated than the shuffle product on $\C(H^{*}(S^{3}\vee\dots\vee S^{2d+1}))$.
\end{ex}
\noindent{\bf Proof of Theorem \ref{modele des lacets libres sur suspension}}
By Lemma \ref{Theoreme de Bott-Samelson} i) and Theorem \ref{Theoreme
  de Goodwillie} (same proof as Theorem \ref{cohomologie lacets libres
  et HAH}), there is a natural DGC quasi-isomorphism
$$\C(\TA\overline{S_*(X)})\buildrel{\simeq}\over\rightarrow S_*((\Sigma
X)^{S^{1}}).$$
Since the cyclic bar construction (Property \ref{cyclic bar
  construction preserve quasi-isomorphisms}), dualizing and
tensorization
preserve quasi-isomorphisms between ${\Bbbk}$-free chain complexes,
$\C(\TA\overline{S_*(X)})^{\vee}$ is weakly equivalent as cochain algebras to
$\C(\TA\overline{C})^{\vee}$. By Lemma \ref{suspension partie
  algebrique}, there is a DGA quasi-isomorphism
$$\C(\TA\overline{C})^{\vee}\buildrel{\simeq}\over\rightarrow
\C({\Bbbk}\oplus s^{-1}\overline{C^{\vee}}).$$
Therefore $\C({\Bbbk}\oplus s^{-1}\overline{C^{\vee}})$ is a weakly DGA
equivalent to $S^{*}((\Sigma X)^{S^{1}})$.\QED

We give now the module structure of the Hochschild homology of
a graded algebra with trivial
product.
Let $V$ be a graded module. The circular permutation toward right
$\tau$ acts on $\T^{n}V$ by
$$\tau.[v_1\vert\dots\vert v_n]=(-1)^{\vert v_n\vert\vert v_1\dots
  v_{n-1}\vert}
[v_n\vert v_1\vert\dots\vert v_{n-1}].$$
Define the {\it invariants} by
$$\T^{n}V^{\tau}=\left\{x\in \T^{n}V,\tau.x=x\right\},
\T V^{\tau}=\oplus_{n=0}^{\infty}\T^{n}V^{\tau}$$ 
and the {\it coinvariants} by
$$\T^{n}V_{\tau}=\frac{\T^{n}V}{\left\{x-\tau.x,x\in \T^{n}V\right\}},
\T V_{\tau}=\oplus_{n=0}^{\infty}\T^{n}V_{\tau}$$
Consider the graded algebra $({\Bbbk}\oplus s^{-1}V)$ with
trivial product.
The Hochschild homology of $({\Bbbk}\oplus s^{-1}V)$
is the sum of the invariants and of the coinvariants of positive
length desuspended:  
$$HH_*( {\Bbbk}\oplus s^{-1}V)=\T V^{\tau}\oplus
 (s^{-1}\otimes 1)\overline{\T V_{\tau}}.$$

Assume now that $V$ is ${\Bbbk}$-free with basis $\beta$.
Then the set of words of length $n$ on $\beta$, denoted $\beta_n$,
is a basis for $\T^{n}V=V^{\otimes n}$.
The circular permutation $\tau$ acts obviously on $\beta_n$ without
sign:
$$\tau*v_1\dots v_n=v_nv_1\dots v_{n-1}.$$
Consider the quotient of the set $\beta_n$ by the cyclic group
generated by $\tau$, $\langle \tau \rangle$.
We denote this quotient set $\langle \tau \rangle\backslash\beta_n $.
For any word $v_1\dots v_n$ denote by $k$ the smallest integer such
that $\tau ^{k}*v_1\dots v_n=v_1\dots v_n$.
Of course the element of $V^{\otimes n}$,
$\tau ^{k}.[v_1|\dots |v_n]$, is either
$+[v_1|\dots |v_n]$ or $-[v_1|\dots |v_n]$.
For any $v_1\dots v_n$ in $\beta_n$, denote by $sym(v_1\dots v_n)$ the
element of $V^{\otimes n}$:
$$\sum_{i=0}^{k-1} \tau ^{i}.[v_1|\dots |v_n].$$
As modules, $\T^{n}V^{\tau}$ is the direct sum
$$
\bigoplus_{\overline{v_1\dots v_n}\in \langle \tau \rangle\backslash\beta_n}
\begin{cases}
{\Bbbk} sym(v_1\dots v_n) &
\text{if $\tau ^{k}.[v_1|\dots |v_n]=+[v_1|\dots |v_n]$},\\
 _2{\Bbbk} sym(v_1\dots v_n) &
\text{if $\tau ^{k}.[v_1|\dots |v_n]=-[v_1|\dots |v_n]$}
\end{cases}
$$
and $\T^{n}V_{\tau}$ is
$$
\bigoplus_{\overline{v_1\dots v_n}\in \langle \tau \rangle\backslash\beta_n}
\begin{cases}
{\Bbbk} \overline{v_1\dots v_n} &
\text{if $\tau ^{k}.[v_1|\dots |v_n]=+[v_1|\dots |v_n]$},\\
 \frac{{\Bbbk}}{2{\Bbbk}} \overline{v_1\dots v_n} &
\text{if $\tau ^{k}.[v_1|\dots |v_n]=-[v_1|\dots |v_n]$}
\end{cases}
$$
In particular, if $V$ is concentrated in even degree or $2=0$ in
${\Bbbk}$ then $HH_*({\Bbbk}\oplus s^{-1}V)$ is ${\Bbbk}$-free.

We suppose now that $V$ is ${\Bbbk}$-free of finite type.
Using Roos direct calculation of the dimension of $\T^{n}V^{\tau}$ when
${\Bbbk}=\mathbb{Q}$ \cite[p. 179-80]{Roos}, we see (Compare \cite[Theorems
1.2.1 and 1.2.2]{Parhizgar}) that the cardinal of
$$\left\{\overline{v_1\dots v_n}\in\langle \tau
  \rangle\backslash\beta_n\mbox{ such that }
\tau ^{k}.[v_1|\dots |v_n]=+[v_1|\dots |v_n]
\mbox{ in }V\otimes_{\Bbbk}\dots\otimes_{\Bbbk} V\right\}$$
is given, when $1\neq-1$ in ${\Bbbk}$, by
$$\frac{1}{n}\sum_{i=1}^{n} \sum_{v_1\dots v_d\in \beta_d}
\varepsilon(\tau^{i},[v_1|\dots |v_d|v_1|\dots |v_d|\cdots |v_1|\dots |v_d]).$$
Here $d$ is the greatest common divisor of $i$ and $n$.
And the integer
$$\varepsilon(\tau^{i},[v_1|\dots |v_d|v_1|\dots |v_d|\cdots |v_1|\dots |v_d])$$
is the sign given by the Koszul rule derived from the action of the
permutation $\tau^{i}$ on the element
of length $n$,
$[v_1|\dots |v_d|v_1|\dots |v_d|\cdots |v_1|\dots |v_d]$.
In particular, by supposing that $V$ is concentrated in even degree,
we see that the cardinal of $\langle \tau \rangle\backslash\beta_n$
is $$\frac{1}{n}\sum_{i=1}^{n} (\mbox{dim }V)^{d}.$$

Let $X$ be a path connected space such that $H_*(X)$ is ${\Bbbk}$-free of finite
type.
The Hopf algebra on $\TC H^{+}(X)$ obtained by tensorization of the
algebra $H^{*}(X)$ is naturally isomorphic as Hopf algebras
to the loop space cohomology $H^{*}(\Omega\Sigma X)$.
The Hochschild homology of $H^{*}(\Sigma X)$,
$HH_*( {\Bbbk}\oplus s^{-1}H^{+}(X))=\T H^{+}(X)^{\tau}\oplus
 (s^{-1}\otimes 1)\overline{\T H^{+}(X)_{\tau}}$ is naturally isomorphic
as graded modules to the free loop space cohomology
$H^{*}\left((\Sigma X)^{S^{1}}\right)$.
This isomorphism of modules is in fact an isomorphism of algebras:
\begin{theor}\label{cohomologie des lacets libres sur
    suspension}
Assume the above hypothesis. The invariants of $\T H^{+}(X)$ form
a graded subalgebra, denoted $\TC H^{+}(X)^{\tau}$,
of the loop space cohomology $H^{*}(\Omega\Sigma X)$.
Consider the
$\left(\TC H^{+}(X)^{\tau},\TC H^{+}(X)^{\tau}\right)$-bimodule
structure on $(s^{-1}\otimes 1)\overline{\T H^{+}(X)_{\tau}}$
induced by the structure of $\left(\TC H^{+}(X),\TC H^{+}(X)\right)$-bimodule
defined on $s^{-1}H^{+}(X)\otimes \T H^{+}(X)$ in
Lemma \ref{suspension partie algebrique}.
Then the associated trivial extension of
$\TC H^{+}(X)^{\tau}$ by $(s^{-1}\otimes
1)\overline{\T H^{+}(X)_{\tau}}$ is naturally isomorphic as graded algebras to the
free loop space cohomology of the suspension of $X$, 
$H^{*}\left((\Sigma X)^{S^{1}}\right)$.
\end{theor}
\pf
Using Lemma \ref{Theoreme de Bott-Samelson} ii),
Theorem \ref{cohomologie lacets libres et HAH}
and Lemma \ref{suspension partie algebrique} with $C=H_*(X)$,
we obtain that the cyclic bar construction on $H^{*}(\Sigma X)$,
$\C({\Bbbk}\oplus s^{-1}H^{+}(X))$ equipped with the product of the
trivial extension given by Lemma \ref{suspension partie algebrique}
for $A=H^{*}(X)$ has the same cohomology algebra as
$H^{*}\left((\Sigma X)^{S^{1}}\right)$.\QED
\begin{rem}
Theorem \ref{cohomologie des lacets libres sur suspension}
claims that the algebra $H^{*}\left((\Sigma X)^{S^{1}}\right)$
depends functorially of the algebra $H^{*}(X)$.
But it is useful to remember that
$H^{*}\left((\Sigma X)^{S^{1}}\right)$ depends functorially of the
Hopf algebra structure of the loop space homology
$H_*(\Omega\Sigma X)=\TA H_+(X)$.

For example, if we return to Example
\ref{Comparaison espace projectif et wedge de spheres},
we obtain the weak equivalences of cochain algebras
$$S^{*}\left((\Sigma \mathbb{CP}^{d})^{S^{1}}\right)\sim
\C\left(\TA H_+(\mathbb{CP}^{d})\right)^{\vee}\sim
\C H^{*}(\Sigma \mathbb{CP}^{d})$$
and
$$S^{*}\left((S^{3}\vee\dots\vee S^{2d+1})^{S^{1}}\right)\sim
\C\left(\TA H_+(S^{2}\vee\dots\vee S^{2d})\right)^{\vee}\sim
\C H^{*}(S^{3}\vee\dots\vee S^{2d+1}).$$
If $\frac{1}{d!}\in {\Bbbk}$ then
$$\TA H_+(\mathbb{CP}^{d})\cong \TA H_+(S^{2}\vee\dots\vee S^{2d})$$
as graded Hopf algebras. We have the isomorphism of cochain algebras
$$\C\left(\TA H_+(\mathbb{CP}^{d})\right)^{\vee}\cong
\C\left(\TA H_+(S^{2}\vee\dots\vee S^{2d})\right)^{\vee}.$$
So finally, when $\frac{1}{d!}\in{\Bbbk}$, we have the
isomorphism of graded Hopf algebras
$$H^{*}(\Omega \Sigma \mathbb{CP}^{d})\cong
H^{*}\left(\Omega (S^{3}\vee\dots\vee S^{2d+1})\right)$$
 and the isomorphism of graded algebras
$$H^{*}\left((\Sigma \mathbb{CP}^{d})^{S^{1}}\right)\cong
H^{*}\left((S^{3}\vee\dots\vee S^{2d+1})^{S^{1}}\right).$$

The converses can be proven easily.
Denote by $x_2$ the generator of $H^{2}(\Omega \Sigma \mathbb{CP}^{d})\cong
H^{2}\left(\Omega (S^{3}\vee\dots\vee S^{2d+1})\right)$.
In $H^{*}(\Omega \Sigma \mathbb{CP}^{d})$, $x_2^{d}\neq0$.
If $d!=0$ in ${\Bbbk}$, $x_2^{d}=0$ in
$H^{*}\left(\Omega (S^{3}\vee\dots\vee S^{2d+1})\right)$.
Therefore, when $d!=0$ in ${\Bbbk}$, there is no isomorphism of graded
algebras between $H^{*}(\Omega \Sigma \mathbb{CP}^{d})$ and
$H^{*}\left(\Omega (S^{3}\vee\dots\vee S^{2d+1})\right)$.
For any space $X$ such that $H_*(X;\mathbb{Z})$ is $\mathbb{Z}$-free of finite
type,
$$H^{*}(X;{\Bbbk})\cong H^{*}(X;\mathbb{Z})\otimes_\mathbb{Z}{\Bbbk}
\mbox{ and so }
H^{*}\left(X;\frac{{\Bbbk}}{d!{\Bbbk}}\right)\cong
H^{*}(X;{\Bbbk})\otimes_{\Bbbk}\frac{{\Bbbk}}{d!{\Bbbk}}$$
as graded algebras.
So we have the implications:
$$H^{*}(\Omega \Sigma \mathbb{CP}^{d};{\Bbbk})\cong
H^{*}\left(\Omega (S^{3}\vee\dots\vee S^{2d+1});{\Bbbk}\right)
\mbox{as graded algebras}$$
$$\Rightarrow
H^{*}(\Omega \Sigma \mathbb{CP}^{d};\frac{{\Bbbk}}{d!{\Bbbk}})\cong
H^{*}\left(\Omega (S^{3}\vee\dots\vee
  S^{2d+1});\frac{{\Bbbk}}{d!{\Bbbk}}\right)\mbox{as graded algebras}$$
$$\Rightarrow
\frac{{\Bbbk}}{d!{\Bbbk}}\mbox{ is the null ring
  }\Rightarrow\frac{1}{d!}\in {\Bbbk}.$$

We prove now that if $d!$ has no inverse in ${\Bbbk}$, there is no
isomorphism of graded algebras between
$H^{*}\left((\Sigma \mathbb{CP}^{d})^{S^{1}}\right)$ and
$H^{*}\left((S^{3}\vee\dots\vee S^{2d+1})^{S^{1}}\right)$.
We have the sequence of isomorphisms of modules \cite[5.3.10]{Loday}
$$H_{*}\left((\Sigma \mathbb{CP}^{d})^{S^{1}}\right)\cong
HH_*\left(\TA(x_2,\dots,x_{2d})\right)\cong
\T(x_2,\dots,x_{2d})_{\tau}\oplus
(s\otimes 1)\overline{\T(x_2,\dots,x_{2d})^{\tau}}.$$
Thus $H_{*}\left((\Sigma \mathbb{CP}^{d})^{S^{1}};\mathbb{Z}\right)$
is $\mathbb{Z}$-free of finite type.
So as for the loop spaces, the proof for the free loop spaces
reduces to the case where $d!=0$ in ${\Bbbk}$.
For $X=\Sigma \mathbb{CP}^{d}$ or $S^{2}\vee\dots\vee S^{2d}$,
by Serre spectral sequence, the inclusion
$\Omega X\hookrightarrow X^{S^{1}}$ induces in cohomology
an isomorphism in degree $2$ and an monomorphism in even degree.
Therefore if $d!=0$ in ${\Bbbk}$, $x_2^{d}\neq0$
in $H^{*}\left((\Sigma \mathbb{CP}^{d})^{S^{1}}\right)$ whereas
$x_2^{d}=0$ in
$H^{*}\left((S^{3}\vee\dots\vee S^{2d+1})^{S^{1}}\right)$.
\end{rem}

It is worth noting the following particular case of
Theorem \ref{cohomologie des lacets libres sur suspension}.
The Hochschild homology of $H^{*}(\Sigma X)$,
$HH_*\left(H^{*}(\Sigma X)\right)$, has a natural structure of
commutative
graded algebra since $H^{*}(\Sigma X)$ is a CDGA.
\begin{cor}\label{suspension et cup produit trivial}
Let $X$ be a path connected space such that $H_*(X)$ is ${\Bbbk}$-free of finite
type. If the cup product on $H^{*}(X)$
is trivial, then $H^{*}\left((\Sigma X)^{S^{1}}\right)$ is naturally isomorphic as
graded algebras to $HH_*\left(H^{*}(\Sigma X)\right)$.
\end{cor}
This Corollary of Theorem
\ref{cohomologie des lacets libres sur suspension}
is proved more easily by applying just Theorem 
\ref{cohomologie lacets libres et HAH},
Lemma \ref{Theoreme de Bott-Samelson} ii) and
Theorem \ref{shc formalite et homologie de Hochschild}.
\section{The Hochschild homology of a commutative
  algebra}\label{exemples Hochschild homology d'une algebre commutative}
If a HAH model of a path connected pointed space $X$ is the cobar
construction on a cocommutative chain coalgebra $C$  
$\Bbbk$-free of finite type such that $C={\Bbbk}\oplus C_{\geq 2}$,
by Theorems \ref{cohomologie lacets libres et HAH} and
\ref{shc formalite et homologie de Hochschild}, the free loop space
cohomology of $X$, $H^{*}(X^{S^{1}})$, is isomorphic as graded algebras
to the Hochschild homology of the CDGA $C^{\vee}$.
In this section, we give various examples of such a space $X$.

Denote by $A$ the cochain algebra $C^{\vee}$.
We suppose now that $A$ is strictly commutative
(i. e. $a^{2}=0$ if $a\in A_{odd}$) and that $\overline{A}$ is
${\Bbbk}$-semifree.
We start by giving a method as general as possible to compute
the Hochschild homology of $A$.

Let $V$ be a graded module.
The free strictly commutative graded algebra on $V$ is denoted
$\Lambda V$.
A {\it decomposable Sullivan Model} of $A$ is a cochain algebra of
the form $(\Lambda V,d)$ where $V=\{V^{i}\}_{i\geq 2}$ is
${\Bbbk}$-free of finite type and $d(V)\subset\Lambda^{\geq 2} V$,
equipped with a quasi-isomorphism of cochains algebras
$(\Lambda V,d)\buildrel{\simeq}\over\rightarrow A$.
If ${\Bbbk}$ is a principal ideal domain, by Theorem 7.1 of
\cite{universal enveloping algebras}, $A$ admits a minimal Sullivan
model.
When ${\Bbbk}$ is a field, minimal Sullivan models are the
decomposable ones \cite[Remark 7.3 i)]{universal enveloping algebras}.

Anyway, suppose now that we have somehow obtained a decomposable
Sullivan model $(\Lambda V,d)$ of $A$ over our arbitrary commutative
ring ${\Bbbk}$. Proposition 1.9 of \cite{radical}
(See also \cite[p 320-2]{Steve et Micheline}) is valid over any
commutative ring ${\Bbbk}$.
Therefore consider the multiplication of $(\Lambda V,d)$:
$$\mu:(\Lambda V',d)\otimes (\Lambda V",d)\rightarrow (\Lambda V,d).$$
By induction on the degree of $V$, we can construct a factorization of
$\mu$:
$$(\Lambda V',d)\otimes(\Lambda V",d)\buildrel{i}\over\rightarrowtail
(\Lambda V'\otimes \Lambda V"
\otimes\Gamma sV,D)\build\twoheadrightarrow_\phi^{\simeq}
(\Lambda V,d)$$
such that
\begin{itemize}
\item[(i)] $D(sv)-(v'-v")\in\Lambda (V^{<n})\otimes \Lambda
  (V^{<n})\otimes\Gamma s(V^{<n})$ for $v\in V^{n}$,
\item[(ii)]$D(\gamma^{k}(sv)=D(sv)\gamma^{k-1}(sv)$ for $v\in
  V^{odd}$ and
\item[(iii)]$\phi(\Gamma(sV)^{+})=0$.
\end{itemize}
Moreover, any such factorization satisfies
\begin{itemize}
\item[(iv)] $i$ is an inclusion of CDGA's such that
$(\Lambda V'\otimes \Lambda V"\otimes\Gamma sV,D)$ is $(\Lambda
V',d)\otimes(\Lambda V",d)$-semifree,
\item[(v)] $\phi$ is a CDGA quasi-isomorphism and
\item[(vi)] $\mbox{Im} D\subset (\Lambda V'\otimes \Lambda V")^{+}\otimes\Gamma sV$.
\end{itemize}

By push out in the category of CDGA's, the multiplication of $A$
extends to a CDGA quasi-isomorphism from
$$(A\otimes A)\otimes_{(\Lambda V,d)\otimes (\Lambda V,d)}(\Lambda V'\otimes \Lambda V"\otimes\Gamma sV,D),$$ which is $A\otimes A$-semifree, to $A$.
The multiplication of $A$ also extends to a CDGA quasi-isomorphism
$\B(A;A;A)\buildrel{\simeq}\over\rightarrow A$.
Since $\overline{A}$ is ${\Bbbk}$-semifree, the bar resolution
$\B(A;A;A)$ is also $A\otimes A$-semifree.
Therefore the cochains algebras
\begin{align*}
\C(A)=&A\otimes_{A\otimes A} \B(A;A;A),\\
(A\otimes\Gamma sV,\overline{D})=&
A\otimes_{A\otimes A}(A\otimes A)\otimes_{(\Lambda V,d)\otimes (\Lambda V,d)}(\Lambda V'\otimes \Lambda V"
\otimes\Gamma sV,D)\\
\intertext{and}
(\Lambda V\otimes\Gamma sV,\overline{D})=&
(\Lambda V,d)\otimes_{(\Lambda V,d)\otimes (\Lambda V,d)}(\Lambda V'\otimes \Lambda V"
\otimes\Gamma sV,D)
\end{align*}
are weakly equivalent as CDGA's (\cite[VIII.2.3]{homology}, for details
\cite[Section 8]{article}).
So finally, we have the isomorphisms of graded algebras
$$HH_*(A)\cong H^{*}(A\otimes\Gamma sV,\overline{D})\cong H^{*}(\Lambda V\otimes\Gamma sV,\overline{D}).$$
\begin{proposition}\label{espace projectif}
The free loop space cohomology on the complex projective space
$\mathbb{CP} ^{n}$, $H^{*}((\mathbb{CP} ^{n})^{S^{1}})$,
is isomorphic as graded algebra to the Hochschild homology of
$H^{*}(\mathbb{CP} ^{n})$.
\end{proposition}
The same result (same proof) holds for the quartenionic projective space
$\mathbb{HP} ^{n}$.
\begin{lem}\label{naturalite du model HAH par rapport aux inclusions de CW-complexes}\cite[8.3 and 8.1g),h)]{anick}\cite[2.1]{modeles d'adams-hilton}
\cite[part 1. of Theorem 6.2]{article}
\begin{itemize}
\item[i)] Let $X$ be a simply connected CW-complex. The Adams-Hilton
model of $X$, denoted ${\cal A}(X)$, can be endowed with a structure
of HAH model for $X$.
\item[ii)] Let $X\hookrightarrow Y$ be an inclusion of simply-connected
CW-complexes.
Consider an Adams-Hilton model of $X$, equipped with a HAH model structure
for $X$, ${\cal A}(X)$. Consider an Adams-Hilton model of $Y$, ${\cal A}(Y)$.
Then there is an structure of HAH model for $Y$ on ${\cal A}(Y)$
that extends the HAH model structure for $X$ of ${\cal A}(X)$.
\end{itemize}
\end{lem}
\noindent{\bf Proof of Proposition \ref{espace projectif}}

By induction on $n$, we suppose that the Adams-Hilton model of
$\mathbb{CP} ^{n-1}$ equipped with its HAH model structure is the
cobar construction $\Omega H_*(\mathbb{CP} ^{n-1})$ equipped with the
shuffle diagonal, denoted $\Delta_s$.
Denote by $\Delta$ the diagonal on
$\Omega H_*(\mathbb{CP} ^{n})=\TA(z_1,\dots, z_{2n-1},d_2)$ obtained
by Lemma~\ref{naturalite du model HAH par rapport aux inclusions de CW-complexes} ii).
This diagonal $\Delta$ is different from the shuffle diagonal $\Delta_s$
only on the top generator $z_{2n-1}$.
So $(\Delta_s-\Delta)z_{2n-1}$ is a cycle. Since
$\Omega \mathbb{CP} ^{n}\thickapprox S^{1}\times\Omega S^{2n-1}$,
for degree reason, it is a boundary.
Therefore, we can construct a derivation homotopy from $\Delta_s$
to $\Delta$.\QED

A decomposable Sullivan model of $\frac{{\Bbbk}[x_2]}{x_2^{n+1}}$
is $(\Lambda(x_2,y_{2n+1}),d)$ with $dy_{2n+1}=x_2^{n+1}$.
Using the general method described above, $HH_*(H^{*}(\mathbb{CP} ^{n}))$
is the cohomology algebra of
$$\left(\frac{{\Bbbk}[x_2]}{x_2^{n+1}}\otimes\Gamma sx_2,sy_{2n+1},
\overline{D}\right)$$ with $\overline{D}sy_{2n+1}=sx(n+1)x^{n}$.
Therefore the graded algebra $H^{*}((\mathbb{CP}^{n})^{S^{1}})$
is the module
\begin{equation*}
\begin{split}
{\Bbbk}\oplus \bigoplus_{1\leq p\leq n,\; i\in\mathbb{N}} {\Bbbk}x^{p}\gamma^{i}(sy)
\oplus \bigoplus_{0\leq p\leq n-1,\; i\in\mathbb{N}} {\Bbbk}x^{p}sx\gamma^{i}(sy)\\
\oplus\bigoplus_{i\in\mathbb{N}}\frac{\Bbbk}{(n+1){\Bbbk}}x^{n}sx\gamma^{i}(sy)
\oplus\bigoplus_{i\in\mathbb{N}^{*}} (_{n+1}{\Bbbk})\gamma^{i}(sy)
\end{split}
\end{equation*}
equipped with the obvious products.
When ${\Bbbk}=\mathbb{Z}$, this is exactly
Proposition 15.33 of~\cite{fibrewise homotopy} (Set $\alpha_i=x\gamma^{i}(sy)$
and $\beta_i=sx\gamma^{i}(sy)$ to make the correspondence).
In particular, if $n+1=0$ in ${\Bbbk}$, we obtain the isomorphism of
graded algebras
$$H^{*}((\mathbb{CP}^{n})^{S^{1}})\cong
H^{*}(\mathbb{CP} ^{n})\times H^{*}((\Omega\mathbb{CP}^{n})).$$

To compute the Hochschild cohomology of an universal enveloping
algebra of a Lie algebra is equivalent as to compute the
Hochschild homology of a commutative algebra:

Consider a differential graded Lie algebra
(in the sense of \cite[1.1(i)]{universal enveloping algebras}) $L$
such that $L=\{L_i\}_{i\geq 1}$ is ${\Bbbk}$-free of finite type.
The universal envelopping algebra of $L$, denoted $UL$,
has a natural structure of DGH~\cite[1.1(ii)]{universal enveloping algebras}.
If $\frac {1}{2}\in {\Bbbk}$, the reduced bar construction
$\B(UL)$ contains a quasi-isomorphic sub-DGC ${\cal C}_*(L)=(\Gamma sL,d_1+d_2)$.
Its dual, denoted ${\mathcal C}^{*}(L)$, is a CDGA called
the {\it (reduced) Cartan-Chevalley-Eilenberg complex}.
Since ${\cal C}_*(L)$ is cocommutative, the cobar construction
$\Omega {\mathcal C}_{*}(L)$ equipped with the shuffle diagonal
is a DGH.
The composite of natural DGA quasi-isomorphisms
$$\Omega {\mathcal C}_{*}(L)\buildrel{\simeq}\over\rightarrow
\Omega\B(UL)\buildrel{\simeq}\over\rightarrow UL$$
is a DGH morphism.
By Theorem~\ref{shc formalite et homologie de Hochschild},
we get immediatly
\begin{lem}\label{algebres enveloppantes et commutatives}
Suppose that $\frac {1}{2}\in {\Bbbk}$.
Let $L$ be a differential graded Lie algebra
such that $L=\{L_i\}_{i\geq 1}$ is ${\Bbbk}$-free of finite type.
Then there is a natural isomorphism of commutative graded algebras
$$HH^{*}(UL)\cong HH_*({\mathcal C}^{*}(L)).$$
\end{lem} 

We give now a large class of spaces who admit the
universal enveloping algebra of a differential graded Lie algebra
as an HAH model:

Let $r\geq 1$ be a fixed integer.
Let $p\geq 2$ be an integer (eventually infinite) such that
$\frac {1}{(p-1)!}\in {\Bbbk}$.
Consider a $r$-connected CW-complex $X$ of finite type and
of dimension $\leq rp$.
We want to compute the free loop space cohomology algebra
$H^{*}(X^{S^{1}})$.
If $p=2$ then, by Freudenthal Suspension Theorem and Proposition 27.5
of~\cite{rational homotopy}, $X$ is the suspension of a co-H space.
And so we have already seen in Section~\ref{free loop space on
  a suspension},
particularly
using Corollary~\ref{suspension et cup produit trivial},   
how to compute its free loop space cohomology algebra
$H^{*}(X^{S^{1}})$.
Therefore, we can suppose that $p\neq 2$.
The Adams-Hilton model of $X$, ${\cal A}(X)$,
is a free chain algebra $(\TA V,d)$ on a ${\Bbbk}$-free graded module $V$
concentrated in degrees between $r$ and $rp-1$, endowed with
a structure of HAH model for $X$.
Therefore, by a deep Theorem of Anick~\cite[5.6]{anick},
there exists a free graded submodule $W\subset\TA V$
such that $d(W)$ embeds into the free graded Lie algebra generated by
$W$, $\mathbb{L}W\subset \TA V$
and such that the DGA morphism
$$U(\mathbb{L}W,d)\buildrel{\cong}\over\rightarrow (\TA V,d)$$
is an HAH isomorphism.
This free differential graded Lie algebra $(\mathbb{L}W,d)$ is the
model $\mathbf{L}(X)$ of Construction 8.4 of~\cite{anick}.
By Lemma~\ref{algebres enveloppantes et commutatives},
we have
\begin{theor}\label{domaine d'Anick}
With the above hypothesis and notations, there is
a natural isomorphism of graded algebras
$$H^{*}(X^{S^{1}})\cong HH_*\left({\mathcal C}^{*}(\mathbf{L}(X))\right).$$
\end{theor}
This Theorem extends the rational case done
by~\cite{Sullivan et  Micheline}
and was the first motivation of this paper after~\cite{article}.

\centerline{\it Universit\'e d'Angers, Facult\'e des Sciences,
2 Boulevard Lavoisier, 49045 Angers, FRANCE}
\centerline{\it e-mail:Luc.Menichi@univ-angers.fr}
\end{document}